%% Plain TeX
%%

\newcount\secno
\newcount\prmno
\newif\ifnotfound
\newif\iffound

\def\namedef#1{\expandafter\def\csname #1\endcsname}
\def\nameuse#1{\csname #1\endcsname}

\long\def\ifundefined#1#2#3{\expandafter\ifx\csname
  #1\endcsname\relax#2\else#3\fi}
\def\hwrite#1#2{{\let\the=0\edef\next{\write#1{#2}}\next}}

% Working with lists
\toksdef\ta=0 \toksdef\tb=2
\long\def\leftappenditem#1\to#2{\ta={\\{#1}}\tb=\expandafter{#2}%
                                \edef#2{\the\ta\the\tb}}
\long\def\rightappenditem#1\to#2{\ta={\\{#1}}\tb=\expandafter{#2}%
                                \edef#2{\the\tb\the\ta}}

\def\lop#1\to#2{\expandafter\lopoff#1\lopoff#1#2}
\long\def\lopoff\\#1#2\lopoff#3#4{\def#4{#1}\def#3{#2}}

\def\ismember#1\of#2{\foundfalse{\let\given=#1%
    \def\\##1{\def\next{##1}%
    \ifx\next\given{\global\foundtrue}\fi}#2}}

% Les commandes
\def\section#1{\vskip1truecm
               \global\def\currenvir{section}
               \global\advance\secno by1\global\prmno=0
               {\bf \number\secno. {#1}}
               \vglue\smallskipamount}

\def\subsection{\global\def\currenvir{subsection}
                \global\advance\prmno by1
                \ind{(\number\secno.\number\prmno) }}
\def\subsec{\global\def\currenvir{subsection}
                \global\advance\prmno by1
                { (\number\secno.\number\prmno)\ }}

\def\proclaim#1{\global\advance\prmno by 1
                {\bf #1 \the\secno.\the\prmno$.-$ }}

\long\def\th#1 \enonce#2\endth{%
   \medbreak\proclaim{#1}{\it #2}\global\def\currenvir{th}\smallskip}

\def\rem#1{\global\advance\prmno by 1
{\it #1}\kern4pt\the\secno.\the\prmno$.-$}

% CROSS-REFERENCES
\def\isinlabellist#1\of#2{\notfoundtrue%
   {\def\given{#1}%
    \def\\##1{\def\next{##1}%
    \lop\next\to\za\lop\next\to\zb%
    \ifx\za\given{\zb\global\notfoundfalse}\fi}#2}%
    \ifnotfound{\immediate\write16%
                 {Warning - [Page \the\pageno] {#1} No reference 
found}}%
                \fi}%
\def\ref#1{\ifx\labellist\empty{\immediate\write16
                 {Warning - No references found at all.}}
               \else{\isinlabellist{#1}\of\labellist}\fi}

\def\newlabel#1#2{\rightappenditem{\\{#1}\\{#2}}\to\labellist}
\def\labellist{}

\def\label#1{%
  \def\given{th}%
  \ifx\given\currenvir%
  {\hwrite\lbl{\string\newlabel{#1}{\number\secno.\number\prmno}}}\fi%
  \def\given{section}%
  \ifx\given\currenvir%
    {\hwrite\lbl{\string\newlabel{#1}{\number\secno}}}\fi%
  \def\given{subsection}%
  \ifx\given\currenvir%
 {\hwrite\lbl{\string\newlabel{#1}{\number\secno.\number\prmno}}}\fi%
  \def\given{subsubsection}%
  \ifx\given\currenvir%
  {\hwrite\lbl{\string%
\newlabel{#1}{\number\secno.\number\subsecno.\number\subsubsecno}}}\fi
  \ignorespaces}
\newwrite\lbl
\newread\testfile
\def\lookatfile#1{\openin\testfile=\jobname.#1
    \ifeof\testfile{\immediate\openout\nameuse{#1}\jobname.#1
                    \write\nameuse{#1}{}
                    \immediate\closeout\nameuse{#1}}\fi%
    \immediate\closein\testfile}%

\def\begin{\newlabel{ACM}{1.1}
\newlabel{min}{1.5}
\newlabel{aut}{1.6}
\newlabel{hyper}{1.7}
\newlabel{smoothdet}{1.8}
\newlabel{cursur}{1.10}
\newlabel{linear}{1.11}
\newlabel{lindet}{1.12}
\newlabel{autsym}{2.3}
\newlabel{symdet}{2.4}
\newlabel{g-1}{3.1}
\newlabel{qtdet}{3.3}
\newlabel{birat}{3.4}
\newlabel{gen}{3.5}
\newlabel{junirat}{3.6}
\newlabel{d/2}{3.7}
\newlabel{theta}{4.2}
\newlabel{simple}{5.3}
\newlabel{AN}{6.1}
\newlabel{surdet}{6.2}
\newlabel{dual}{6.3}
\newlabel{minors}{6.7}
\newlabel{spf}{7.1}
\newlabel{surpf}{7.2}
\newlabel{charp}{7.3}
\newlabel{expf}{7.4}
\newlabel{gensurpf}{7.6}
\newlabel{n-folds}{8.2}
\newlabel{exist}{8.3}
\newlabel{ijunirat}{8.8}
\newlabel{3f}{8.9}}

%% Fin de la numerotation automatique

\magnification 1250
\pretolerance=500 \tolerance=1000  \brokenpenalty=5000
\mathcode`A="7041 \mathcode`B="7042 \mathcode`C="7043
\mathcode`D="7044 \mathcode`E="7045 \mathcode`F="7046
\mathcode`G="7047 \mathcode`H="7048 \mathcode`I="7049
\mathcode`J="704A \mathcode`K="704B \mathcode`L="704C
\mathcode`M="704D \mathcode`N="704E \mathcode`O="704F
\mathcode`P="7050 \mathcode`Q="7051 \mathcode`R="7052
\mathcode`S="7053 \mathcode`T="7054 \mathcode`U="7055
\mathcode`V="7056 \mathcode`W="7057 \mathcode`X="7058
\mathcode`Y="7059 \mathcode`Z="705A
\def\spacedmath#1{\def\packedmath##1${\bgroup\mathsurround =0pt##1
\egroup$}
\mathsurround#1
\everymath={\packedmath}\everydisplay={\mathsurround=0pt}}
\def\nospacedmath{\mathsurround=0pt
\everymath={}\everydisplay={} } \spacedmath{2pt}
\def\qfl#1{\buildrel {#1}\over {\longrightarrow}}
\def\phfl#1#2{\normalbaselines{\baselineskip=0pt
\lineskip=10truept\lineskiplimit=1truept}\nospacedmath\smash{\mathop
{\hbox to 8truemm{\rightarrowfill}}
\limits^{\scriptstyle#1}_{\scriptstyle#2}}}
\def\hfl#1#2{\normalbaselines{\baselineskip=0truept
\lineskip=10truept\lineskiplimit=1truept}\nospacedmath\smash{\mathop
{\hbox to
12truemm{\rightarrowfill}}\limits^{\scriptstyle#1}_{\scriptstyle#2}}}
\def\diagram#1{\def\normalbaselines{\baselineskip=0truept
\lineskip=10truept\lineskiplimit=1truept}   \matrix{#1}}
\def\vfl#1#2{\llap{$\scriptstyle#1$}\left\downarrow\vbox to
6truemm{}\right.\rlap{$\scriptstyle#2$}}
\font\eightrm=cmr8         \font\eighti=cmmi8
\font\eightsy=cmsy8        \font\eightbf=cmbx8
\font\eighttt=cmtt8        \font\eightit=cmti8
\font\eightsl=cmsl8        \font\sixrm=cmr6
\font\sixi=cmmi6           \font\sixsy=cmsy6
\font\sixbf=cmbx6\catcode`\@=11
\def\eightpoint{%
  \textfont0=\eightrm \scriptfont0=\sixrm \scriptscriptfont0=\fiverm
\def\rm{\fam\z@\eightrm}%
  \textfont1=\eighti  \scriptfont1=\sixi 
\scriptscriptfont1=\fivei
  \textfont2=\eightsy \scriptfont2=\sixsy \scriptscriptfont2=\fivesy
  \textfont\itfam=\eightit
  \def\it{\fam\itfam\eightit}%
  \textfont\slfam=\eightsl
  \def\sl{\fam\slfam\eightsl}%
  \textfont\bffam=\eightbf \scriptfont\bffam=\sixbf
  \scriptscriptfont\bffam=\fivebf
  \def\bf{\fam\bffam\eightbf}%
  \textfont\ttfam=\eighttt
  \def\tt{\fam\ttfam\eighttt}%
  \abovedisplayskip=9pt plus 3pt minus 9pt
  \belowdisplayskip=\abovedisplayskip
  \abovedisplayshortskip=0pt plus 3pt
  \belowdisplayshortskip=3pt plus 3pt 
  \smallskipamount=2pt plus 1pt minus 1pt
  \medskipamount=4pt plus 2pt minus 1pt
  \bigskipamount=9pt plus 3pt minus 3pt
  \normalbaselineskip=9pt
  \normalbaselines\rm}\catcode`\@=12 
\newcount\noteno
\noteno=0 
\def\up#1{\raise 1ex\hbox{\sevenrm#1}}
\def\note#1{\global\advance\noteno by1
\footnote{\parindent0.4cm\kern2pt\up{\number\noteno}\
}{\vtop{\eightpoint\baselineskip12pt\hsize15.5truecm\noindent
#1}}\parindent 0cm}
\def\mono{\lhook\joinrel\mathrel{\longrightarrow}}
\def\sdir_#1^#2{\mathrel{\mathop{\kern0pt\oplus}\limits_{#1}^{#2}}}
\def\pprod_#1^#2{\raise2pt
\hbox{$\mathrel{\scriptstyle\mathop{\kern0pt\prod}\limits_{#1}^{#2}}$}}
\def\pc#1{\tenrm#1\sevenrm}
\def\tx{\kern-1.5pt -}
\def\cqfd{\kern 2truemm\unskip\penalty 500\vrule height 4pt depth 0pt 
width 4pt} 
\def\virg{\raise
.4ex\hbox{,}}
\def\decale#1{\smallbreak\hskip 28pt\llap{#1}\kern 5pt}
\def\no{n\up{o}\kern 2pt}
\def\ind{\par\hskip 1truecm\relax}
\def\indp{\par\hskip 0.5truecm\relax}
\def\moins{\mathrel{\hbox{\vrule height 3pt depth -2pt width 6pt}}}
\def\rond{\kern 1pt{\scriptstyle\circ}\kern 1pt}
\def\prond{\kern .3pt{\scriptscriptstyle\circ}\kern .6pt}
\def\iso{\vbox{\hbox to .8cm{\hfill{$\scriptstyle\sim$}\hfill}
\nointerlineskip\hbox to .8cm{{\hfill$\longrightarrow $\hfill}}
}}
\def\End{\mathop{\rm End}\nolimits}
\def\Hom{\mathop{\rm Hom}\nolimits}
\def\Aut{\mathop{\rm Aut}\nolimits}
\def\im{\mathop{\rm Im}\nolimits}
\def\Ker{\mathop{\rm Ker}\nolimits}
\def\Coker{\mathop{\rm Coker}}
\def\det{\mathop{\rm det}\nolimits}
\def\Pic{\mathop{\rm Pic}\nolimits}

\def\dim{\mathop{\rm dim}\nolimits}
\def\Card{\mathop{\rm Card}\nolimits}

\def\rk{\mathop{\rm rk\,}\nolimits}

\def\ch{\mathop{\rm char}\nolimits}
\def\pf{\mathop{\rm pf}\nolimits}
\def\ten{^{\scriptscriptstyle \otimes 2}}
\def\pr{{\it Proof}\kern2.3pt: }
\font\san=cmssdc10
\def\ext{\hbox{\san \char3}}
\def\s{\hbox{\san \char83}}
\def\id{\hbox{\san \char73}}
\def\r{\hbox{\san \char82}}

\def\pn{{\bf P}^n}
\def\pl{{\bf P}^2}
\def\pt{{\bf P}^3}

\font\ru=wncyb10
\def\gf{\hbox{\ru G}_*({\cal F})}
\font\scr=eusm10\font\sbr=eurm10
\def\en{\hbox{\scr\char69\hskip-0.5pt\sbr\char110\hskip-0.5pt\char100}}

\input amssym.def
\input amssym
\vsize = 25.4truecm
\hsize = 16truecm
%\hoffset = -.15truecm
\voffset = -.7truecm
\parindent=0cm
\baselineskip15pt
\overfullrule=0pt
\frenchspacing 
\begin
\centerline{\bf  Determinantal hypersurfaces}
\smallskip
\smallskip \centerline{Arnaud {\pc BEAUVILLE}} 
\vskip0.8cm
{\leftskip=12truecm \it To Bill Fulton\par}\bigskip
{\bf Introduction}
\smallskip
\subsection We discuss in this paper  which homogeneous form on
$\pn$ can be written as the determinant of a matrix with homogeneous
entries (possibly  symmetric), or the
pfaffian of  a skew-symmetric matrix. This question has been
considered in various particular cases (see the historical comments
below), and we believe that the general result is well-known from
the experts; but we have been unable to find it in
the literature. The aim of this paper is to fill this gap.
\ind We will discuss at the outset the general structure
theorems; roughly, they show that expressing a homogeneous form
$F$ as a determinant (resp. a pfaffian) is equivalent to produce
a line bundle (resp. a rank 2 vector bundle) of a certain type on
the hypersurface $F=0$. The rest of the paper consists of
applications. We have restricted our attention to  {\it smooth}
hypersurfaces; in fact we are particularly interested in the case
when the {\it generic} form of degree $d$ in $\pn$ can be written
in one of the above forms.  When this is the case,  the moduli
space of pairs $(X,E)$, where $X$ is a smooth hypersurface of
degree $d$ in $\pn$ and $E$ a rank 1 or 2 vector bundle
satisfying appropriate conditions, appears as a quotient of an
open subset of a certain vector space of matrices; in particular,
this moduli space is {\it unirational}. This is the case for
instance of the universal family  of Jacobians   of plane curves
(Cor. \ref{junirat}), or  of intermediate Jacobians of cubic
threefolds  (Cor. \ref{ijunirat}).
\ind Unfortunately  this situation does not occur too frequently: we 
will show that only  curves and cubic surfaces admit generically a
determinantal equation. The situation is slightly better for
pfaffians: plane curves of any degree,  surfaces of degree $\le
15$ and  threefolds of degree $\le 5$ can be generically defined by a
linear pfaffian.
\medskip
\subsec{\it Historical comments}\ind The representation of curves and surfaces
of small degree as
linear determinants is a  classical subject. The case of
cubic surfaces was already known in the middle of the last
century  [G]; other examples of curves and surfaces are treated
in [S]. The  general homogeneous forms which can
be expressed as linear  determinants are determined in
[D]. A modern treatment for  plane curves appears in [C-T];  the
result has been rediscovered a number of times since then.
	\ind The representation of the
plane quartic as a symmetric determinant goes back again to
1855 [H]; plane curves of any degree are treated in [Di]. Cubic
and quartic surfaces defined by linear symmetric determinants
(``symmetroids")  have been also studied early  [Ca]. Surfaces of
 arbitrary degree are thoroughly treated in [C1]; an overview of the
use of symmetric resolutions can be found in [C2].
\ind Finally, the only reference we know about pfaffians is Adler's
proof  that a generic cubic threefold can be written as a linear
pfaffian ([A-R], App. V).

 \medskip
\subsec{\it Conventions}
\ind We work over an arbitrary field $k$, not necessarily
algebraically closed. Unless explicitely stated, all geometric objects
are defined over $k$.
\vskip0.4truecm
\ind  {\eightpoint{\it Acknowledgements}: I thank F. Catanese for his useful
comments, and F.-O. Schreyer for providing the computer-aided
proof of Prop. \ref{gensurpf} b) and 8.9 below (see
Appendix).}
\section{General results: determinants}
\subsection\label{ACM} Let 
 ${\cal F}$ be a coherent sheaf on $\pn$. We will say that ${\cal
F}$ is {\it arithmetically Cohen-Macaulay} (ACM for short) if:
\indp a) ${\cal F}$ is Cohen-Macaulay, that is, the ${\cal O}_x$\tx 
module ${\cal F}_x$ is Cohen-Macaulay for every $x$ in $\pn$;
\indp b) $H^i(\pn,{\cal F}(j))=0$ for $1\le i\le \dim({\rm
Supp}\,{\cal F})-1$ and 
$j\in {\bf Z}$.
\ind Put $\s^n=k[X_0,\ldots,X_n]=\sdir_{j\in{\bf Z}}^{}H^0(\pn,
{\cal O}_{\pn}(j))$ (we will often drop the superscript $n$ if there
 is no danger of confusion). Following EGA, we denote by
$\gf$ the $\s$\tx module
$\sdir_{j\in {\bf Z}}^{}H^0(\pn, {\cal F}(j))$. The following
well-known remark explains the terminology:
\th Proposition
\enonce The sheaf ${\cal F}$ is {\rm ACM} if and only if the $\s$\tx 
module $\gf$ is Cohen-Macaulay.
\endth
\pr  Let $U:={\bf A}^{n+1}\moins\{0\}$.
The projection $p:U\rightarrow \pn$ is affine, and
satisfies $p_*{\cal O}_U=\sdir_{j\in{\bf Z}}^{}{\cal O}_{\pn}(j)$.
The $\s$\tx module $\gf$ defines a coherent
sheaf $\widetilde{\cal F}$ on ${\bf A}^{n+1}$, whose restriction to $U$ 
is isomorphic to $p^*{\cal F}$. Therefore  $H^i(U,\widetilde{\cal F})$
is isomorphic to
$\sdir_{j\in{\bf Z}}^{}H^i(\pn,{\cal F}(j)$. The long exact sequence 
of local cohomology
$$\cdots\longrightarrow  H^i_{\{0\}}({\bf A}^{n+1},\widetilde{\cal F})
\longrightarrow H^i({\bf A}^{n+1},\widetilde{\cal F}) \longrightarrow
H^i(U,\widetilde{\cal F})\longrightarrow \cdots$$
implies $H^0_{\{0\}}({\bf A}^{n+1},\widetilde{\cal F})=H^1_{\{0\}}({\bf
A}^{n+1},\widetilde{\cal F})=0$, and give isomorphisms
 $$\sdir_{j\in {\bf Z}}^{}H^i(\pn,{\cal F}(j))\iso H^{i+1}_{\{0\}}({\bf
A}^{n+1},\widetilde{\cal F})\quad{\rm for}\  i\ge 1\ .$$
Thus condition b) of (\ref{ACM}) is equivalent to
$H^i_{\{0\}}(\widetilde{\cal F})=0$ for $i<\dim(\widetilde{\cal F})$,
that is to $\widetilde{\cal F}_0$ being Cohen-Macaulay. On the other
hand, since $p$ is smooth, condition a) is equivalent to
$\widetilde{\cal F}_v$ being Cohen-Macaulay for all $v\in U$, hence
the Proposition.\cqfd\medskip
\ind Let us mention incidentally the following corollary, due to
Horrocks:
\th Corollary
\enonce A locally free sheaf ${\cal F}$ on $\pn$ with $H^i(\pn,{\cal
F}(j))=0$ for $1\le i\le n-1$ and $j\in{\bf Z}$ splits as a direct sum
of line bundles.
\endth
\pr  The $\s$\tx module $\gf$ is Cohen-Macaulay of maximal
dimension, hence projective; it is therefore free as a $\s$\tx
graded module, that is isomorphic to a direct sum
$\s(d_1)\oplus\cdots\oplus \s(d_r)$ ([Bo], \S 8, Prop. 8). Since
${\cal F}$ is the sheaf on ${\rm Proj}(\s)$ associated to $\gf$, it is 
isomorphic to ${\cal O}_{\pn}(d_1)\oplus\cdots\oplus {\cal
O}_{\pn}(d_r)$.   
\bigbreak
{\bf Theorem A}$.-$ {\it  Let ${\cal F}$ be an {\rm ACM} sheaf on $\pn$,
of dimension $n-1$. There exists an exact sequence
$$0\rightarrow  \sdir_{i=1}^{\ell }{\cal O}_{\pn}(e_i)\qfl{M}
\sdir_{i=1}^{\ell }{\cal O}_{\pn}(d_i)\longrightarrow {\cal
F}\rightarrow 0\ .\eqno(A1)$$
\ind Conversely, let $M: \sdir_{i=1}^{\ell }{\cal O}_{\pn}(e_i)
\rightarrow \sdir_{i=1}^{\ell }{\cal O}_{\pn}(d_i)$ be an injective
homomorphism; the  cokernel of $M$ is {\rm ACM} and  its support  is the
hypersurface $\det M=0$.  
\smallskip 
Proof}:  Suppose that ${\cal F}$ is ACM of dimension $n-1$. The
Cohen-Macaulay $\s$\tx module $\gf$ has projective dimension $1$; by 
Hilbert's theorem ([Bo], \S 8, Cor. 3 of Prop.~8), it admits a free 
graded resolution of the form
$$0\rightarrow \sdir_{i=1}^{\ell }\s(e_i)\longrightarrow 
\sdir_{i=1}^{\ell }\s(d_i)\longrightarrow \gf \rightarrow 0\ 
,\eqno(A2)$$ which gives  $(A1)$ by taking the associated sheaves on
$\pn$.
\ind Conversely, suppose given the exact sequence $(A1)$. The support  
of ${\cal F}$ consists of the points $x$ of $\pn$ where $M(x)$ is not
injective, that is where $\det M(x)=0$. Since $M$ is generically 
injective this is a hypersurface in
$\pn$. 
\ind For every $x\in\pn$,
the ${\cal O}_{\pn,x}$\tx module
${\cal F}_x$ has projective dimension $\le 1$, hence depth $\ge \dim 
{\cal O}_{\pn,x} -1=\dim {\cal F}_x$; thus it is Cohen-Macaulay.
From $(A1)$ we deduce $H^i(\pn,{\cal F}(j))=0$ for $1\le i\le n-2$, 
hence 
${\cal F}$ is ACM.\cqfd\medskip
\subsection The homomorphism $M$ is given by a matrix 
$(m_{ij})\in {\bf M}_\ell (\s)$, with $m_{ij}$ homogeneous of degree
$(d_i-e_j)$; we will  use the same letter $M$ to denote this
matrix.\medskip
\subsection\label{min} Let  ${\cal F}$ be  an ACM sheaf on $\pn$ of
dimension
$n-1$. We will always take for $(A2)$ a {\it minimal} graded free
resolution of
$\gf$: this means that the   images in $\gf$ of the generators of
$\s(d_i)$ $(1\le i\le \ell )$ form a minimal system of generators of the
$\s$\tx module $\gf$. Such a resolution is unique up to isomorphism.
The resolution $(A2)$ is minimal if and only if the matrix $(m_{ij})$
is zero modulo $(X_0,\ldots,X_n)$, that is, if and only if
$m_{ij}=0$ whenever $d_i=e_j$. 
\ind We will refer to the corresponding exact sequence $(A1)$,
slightly abusively, as the {\it minimal resolution} of the sheaf
${\cal F}$.
\smallskip 
\subsection\label{aut} The minimal resolution $0\rightarrow
L_1\rightarrow L_0\rightarrow {\cal F}\rightarrow 0$, with
$L_1=\sdir_{i=1}^{\ell }{\cal O}_{\pn}(e_i)$ and  $L_0=
\sdir_{i=1}^{\ell }{\cal O}_{\pn}(d_i)$, is unique
up to isomorphism, but this isomorphism is not unique in general; {\it
it is unique if}
$\max(e_j)<\min(d_i)$ (in particular in the linear case). Indeed this
condition implies
$\Hom(L_0,L_1)=0$, and therefore the map 
$\End(L_0)\rightarrow \Hom(L_0,{\cal F})$ is injective; thus the only
automorphism of $L_0$ which induces the identity on ${\cal F}$ is the
identity. If moreover every automorphism of ${\cal F}$ is scalar, we
see that the only pairs of automorphisms $P\in \Aut(L_0)$, $Q\in
\Aut(L_1)$ such that $PM=MQ$ are the pairs $(\lambda,\lambda)$ for
$\lambda\in k^*$.
 \medskip
\subsection\label{hyper} In this paper we will mainly use Theorem A in
the following way: we will start from an integral
(usually smooth) hypersurface $X$ and a vector bundle $E$
 of rank $r$ on $X$; we will still say that $E$ is ACM 
if it is so as an ${\cal O}_{\pn}$\tx module, that is,
$H^i(X,{\cal F}(j))=0$ for $1\le i\le n-2$ and $j\in{\bf Z}$. For
such a sheaf Theorem A provides a minimal resolution $(A1)$;
localizing at the generic point of $X$
and using the structure theorem for  matrices over a principal ring we
get $\det M=F^r$, where $F=0$ is an equation of $X$. This gives the
following corollary:
\th Corollary
\enonce Let $X$ be a smooth hypersurface  in $\pn$,
given by an equation $F=0$. 
\ind {\rm a)} Let $L$ be a line bundle on $X$ with
$H^i(X,L(j) )=0$ for $1\le i\le n-2$ and all $j\in{\bf Z}$.
Then $L$ admits a minimal resolution
$$0\rightarrow  \sdir_{i=1}^{\ell }{\cal O}_{\pn}(e_i)\qfl{M} 
\sdir_{i=1}^{\ell }{\cal O}_{\pn}(d_i)\longrightarrow L\rightarrow
0$$with $F=\det M$.
\ind {\rm b)} Conversely, let $M=(m_{ij})\in {\bf M}_\ell (\s)$,  with
$m_{ij}$ homogeneous of degree $(d_i-e_j)$ and $F=\det M$. Then the
cokernel of
$M: \sdir_{i=1}^{\ell }{\cal O}_{\pn}(e_i)\longrightarrow  
\sdir_{i=1}^{\ell }{\cal O}_{\pn}(d_i)$ is a line bundle
$L$ on $X$ with the above properties.\cqfd
\endth\label{smoothdet}
\subsection The apparent generality of this Corollary
is somewhat misleading: taking for $L$ the line bundle ${\cal
O}_X(j)$ gives rise to the trivial case $\ell =1$, $M=(F)$. Thus
 if $\Pic(X)$ is generated by ${\cal O}_X(1)$ the
hypersurface can {\it not} be defined by a $\ell \times \ell $ 
determinant with $\ell >1$. So interesting situations occur only for 
curves and surfaces. In particular, we infer from the
Noether-Lefschetz  theorem that the generic hypersurface of degree
$d$ in ${\bf P}^n$ can be expressed in a non-trivial way as a
determinant only  if $n=2$ or $n=3$ and $d\le 3$. On the other hand we
will
 see in (3.1) and (6.4) that any smooth plane curve or  cubic
surface can be defined by a linear determinant.
\subsection Conversely, given integers $e_i,d_j$, one may ask
whether a general matrix $(m_{ij})\in {\bf M}_\ell (\s)$ with
$\deg m_{ij}=d_i-e_j$ defines a smooth curve or 
 surface. If we order the numbers $e_i,d_j$ so that $e_1\le \ldots
\le e_\ell $ and $d_1\le \ldots \le d_\ell $, {\it a sufficient 
condition is the inequality $d_{i}>e_{i+1}$ for} $1\le i <\ell $.  
Indeed we can consider the matrix $$M=\pmatrix{
F_1   &   G_1&  0&  \cdots & 0\cr
0   &   F_2    & G_2& \cdots & 0\cr
\vdots &  &\ddots &\ddots&\vdots\cr
0 &  & &  F_{\ell -1}& G_{\ell -1}\cr
G_\ell   &0& \cdots &0 &F_\ell  \cr
}$$where the entries are product of linear
forms. Then $\det M$ can be written in the form $\pprod_{}^{}
L_i+\pprod_{}^{} P_i$, where 
$L_i,P_j$ are arbitrary linear forms. We obtain in this way, for
instance, the Fermat hypersurface\note{If  $\ch(k)\mid d$ consider the
surface $X_0^{}(X_0^{d-1}+X_1^{d-1})+
(X_1^{}+X_2^{})(X_2^{d-1}+X_3^{d-1})=0$.} $\sum X_i^d=0$ in $\pl$ or 
$\pt$.
\label{cursur}
\bigskip 
\ind  The integers $e_i,d_j$ which occur in the minimal resolution are
determined by the $\s$\tx module $\gf$; we will see some examples in
the next sections. We will be particularly interested by the case 
when the entries $(m_{ij})$ are linear forms; in this case we will
say for short that the matrix $M$ is {\it linear}. There is a handy
characterization  of the sheaves which give rise to linear matrices:
\th Proposition
\enonce Let ${\cal F}$ be a coherent sheaf on $\pn$. The following 
conditions are equivalent:
\indp{\rm (i)} There exists an exact sequence
$$0\rightarrow {\cal O}_{\pn}(-1)^\ell \longrightarrow {\cal
O}_{\pn}^\ell \longrightarrow {\cal F}\rightarrow 0\ ;$$
\indp{\rm (ii)} ${\cal F}$ is {\rm ACM} of dimension $n-1$, and
$$H^0(\pn,{\cal F}(-1))=H^{n-1}(\pn, {\cal F}(1-n))=0\ .$$
\endth\label{linear}
{\it Proof} : In view of Theorem A the implication ${\rm (i)}
\Rightarrow {\rm (ii)}$ is clear. Assume that (ii) holds; then
$H^i(\pn, {\cal F}(-i))=0$ for  $i\ge 1$, that is, ${\cal F}$ is 
$0$\tx {\it regular} in the sense of Mumford ([Md], lect. 14).
By {\it loc.\ cit}., this implies that ${\cal F}$ is spanned by
its global sections and that  the natural map $$H^0(\pn, {\cal
F}(j))\otimes H^0(\pn, {\cal O}_{\pn}(1))\rightarrow H^0(\pn,
{\cal F}(j+1))$$ is surjective for
$j\ge 0$. Since $H^0(\pn,{\cal F}(-1))=0$, this means that the
multiplication map $\s\otimes_kH^0(\pn,{\cal F})\rightarrow \gf$
is surjective, and therefore  the minimal resolution of ${\cal
F}$ takes the form:
$$0\rightarrow  \sdir_{i=1}^{\ell }{\cal O}_{\pn}(e_i)\qfl{M} 
{\cal O}_{\pn}^\ell \qfl{p} {\cal F}\rightarrow 0$$with $\ell =\dim
H^0(\pn,{\cal F})$. Since $H^0(p)$ is bijective and $H^{n-1}(\pn,
{\cal F}(1-n))=0$, we must have $e_i=-1$ for all $i$.\cqfd\medskip

\ind We can again reformulate this result as: 
\th Corollary
\enonce Let $X$ be a smooth hypersurface of degree $d$ in $\pn$,
given by an equation $F=0$. 
\ind {\rm a)} Let $L$ be a line bundle on $X$ with $H^i(X,L(j) )=0$
for $1\le i\le n-2$ and all $j\in{\bf Z}$, and
$H^{0}(X,L(-1))=H^{n-1}(X,L(1-n) )=0$. There exists a $d\times d$ linear
matrix $M$ such that $F=\det M$, and an exact
sequence
$$0\rightarrow {\cal O}_{\pn}(-1)^d\qfl{M} {\cal
O}_{\pn}^d\longrightarrow L\rightarrow 0\ .$$
\ind {\rm b)} Conversely, let $M$ be a $d\times d$ linear matrix such  
that
$F=\det M$. Then the cokernel of $M: {\cal O}_{\pn}(-1)^d\rightarrow 
{\cal O}_{\pn}^d$ is a line bundle $L$ on
$X$ with the above properties.\cqfd
\endth\label{lindet}
\section{General results: symmetric determinants and pfaffians}
\subsection We will now put an extra data on our ACM sheaf. Let  ${\cal
F}$ be a torsion-free sheaf on an  integral variety $X$, and $L$ a
line bundle on $X$;  a bilinear form $\varphi:{\cal
F}\otimes_{{\cal O}_X}{\cal F}\rightarrow L$ is said to be {\it
invertible} if the associated homomorphism $\kappa:{\cal F}\rightarrow
{\cal H}om_{{\cal O}_X} ({\cal F},L)$ is an isomorphism. We will
consider forms which are $\varepsilon$\tx symmetric $(\varepsilon=\pm
1)$, that is, such that ${}^t\kappa=\varepsilon\,\kappa$. 
\medbreak
{\bf Theorem B}$.-$ {\it Assume $\ch(k)\not= 2$. Let $X$ be
an integral hypersurface of degree
$d$ in $\pn$, and ${\cal F}$ a torsion-free {\rm ACM} sheaf on
$X$, equipped  with an $\varepsilon$\tx
symmetric invertible form
${\cal F}\otimes{\cal F}\rightarrow {\cal O}_X(d+t)$ $(t\in{\bf Z})$.
Then
${\cal F}$ admits a resolution 
$$0\rightarrow L_0^*(t)\qfl{M} L_0\longrightarrow {\cal
F}\rightarrow 0\ ,\eqno(B1)$$
where $L_0=\oplus {\cal O}_{\pn}(d_i)$ and $M$ is $\varepsilon$\tx
symmetric, that is,
${}^tM=\varepsilon\,M$.
\ind Conversely, if a sheaf ${\cal F}$ on $X$ fits into the exact
sequence $(B1)$, it is {\rm ACM}, torsion-free, and admits an 
$\varepsilon$\tx symmetric invertible form}
${\cal F}\otimes{\cal F}\rightarrow {\cal O}_X(d+t)$.
\smallskip 
\pr  Consider a minimal resolution
$$0\rightarrow L_1\qfl{M} L_0 \qfl{p} {\cal F}\rightarrow 0$$ 
of ${\cal F}$. Applying the functor ${\cal H}om^{}_{{\cal
O}_{\pn}}(\,*\, ,{\cal O}_{\pn}(t))$ gives an exact sequence  
$$0\rightarrow L_0^*(t)\qfl{{}^tM} L_1^*(t) \qfl{} {\cal
E}xt^1_{{\cal O}_{\pn}}({\cal F},{\cal O}_X(t))\rightarrow 0$$
and the vanishing of ${\cal E}xt^i_{{\cal O}_{\pn}}({\cal
F},{\cal O}_X(t))$ for $i\not= 1$.
\ind Let $i$ be the  embedding of $X$ into $\pn$; put ${\cal F}'=
{\cal H}om_{{\cal O}_X}({\cal F},{\cal O}_X(d+t))$.
Gro\-then\-dieck duality provides a canonical isomorphism ${\cal
E}xt^1_{{\cal O}_{\pn}}({\cal F},{\cal O}_X(t))\iso
i_*{\cal
 F}'$. Thus the above exact
sequence  gives a minimal resolution of the \hbox{${\cal
O}_{\pn}$\tx module} ${\cal F}'$; the isomorphism $\kappa:{\cal
F}\rightarrow {\cal F}'$ extends to an isomorphism of resolutions: 
$$\diagram{0\rightarrow & L_1 & \hfl{M}{} & L_0 & \hfl{p}{} &{\cal
F} & \rightarrow 0&\cr
& \vfl{B}{} & & \vfl{A}{} & &\vfl{\kappa}{} &&\cr
0\rightarrow & L_0^*(t) &\hfl{{}^tM}{} & L_1^*(t) &
\hfl{q}{} & {\cal F}'
&\rightarrow 0&.\cr }$$
Applying the functor ${\cal H}om^{}_{{\cal
O}_{\pn}}(\,*\, ,{\cal O}_{\pn}(t))$ leads to another commutative
diagram:
$$\diagram{0\rightarrow & L_1 & \hfl{M}{} & L_0 & \hfl{p}{} &{\cal
F} & \rightarrow 0&\cr
& \vfl{{}^t\!A}{} & & \vfl{{}^tB}{} & &\vfl{{}^t\kappa}{} &&\cr
0\rightarrow & L_0^*(t) &\hfl{{}^tM}{} & L_1^*(t) &
\hfl{q}{} & {\cal F}'
&\rightarrow 0&.\cr }$$
Since ${}^t\kappa=\varepsilon\,\kappa$, we have
$q\rond{}^tB= {}^t\kappa\rond p=\varepsilon\,q\rond A$, which means
that there exists a map $C:L_0\rightarrow L_0^*(t)$ such that
${}^tB-\varepsilon\,A={}^tM C$.
\ind Since ${}^tB M ={}^tM \,^t\!A$, we have $${}^tM
C M =({}^tB-\varepsilon\,A) M= {}^t(A
M)-\varepsilon(A M)=-\varepsilon\,{}^tM {}^tC M$$ 
and therefore the map $A':=A+{\varepsilon\over 2}{}^tM C$ satisfies 
${}^t(A' M)=\varepsilon\,A' M$. Moreover we still have
$q\rond A'=\kappa\rond p$, so $A'$ is an isomorphism. We have an
exact sequence
$$0\rightarrow L_0^*(t)\qfl{M'} L_0\qfl{p}{\cal F}\rightarrow 0\
,$$ where $M':=A'^{-1}\,^tM$ satisfies
${}^tM'=\varepsilon\,M'$.
\ind Conversely, starting from the exact sequence $(B1)$, 
Grothendieck duality implies as above an isomorphism
$\kappa:{\cal F}\rightarrow {\cal H}om({\cal F},{\cal O}_X(d+t))$;
applying again the functor ${\cal H}om^{}_{{\cal
O}_{\pn}}(\,*\, ,{\cal O}_{\pn}(t))$ we obtain
${}^t\kappa=\varepsilon\,\kappa$.\cqfd\medskip

\rem{Remark} The result remains valid in characteristic $2$ under the
extra hypo\-thesis $\max(e_j)<\min(d_i)$: indeed, with the above
notation, the relation $q\rond A= q\rond{}^tB$ 
implies then directly $A={}^tB$ (\ref{aut}), and we can take
$M'=A^{\!-1}
 \,^tM$.
\ind F. Catanese pointed out that his proof in [C1] for
symmetric surfaces extends readily to the case considered here; it
has the advantage of working equally well in characteristic $2$,
without the above restriction on the degrees.
\medskip 
\subsection\label{autsym} Assume again $\max(e_j)<\min(d_i)$. Let  
$0\rightarrow P_0^*(t')\qfl{M'} P_0\qfl{p'}{\cal F}\rightarrow 0$
be another resolution (B1) of ${\cal F}$; then  we
have $t=t'$ and a commutative diagram
$$\diagram{0\rightarrow & L_0^*(t) & \hfl{M}{} & L_0 & \hfl{p}{} &{\cal
F} & \rightarrow 0&\cr
& \vfl{B}{} & & \vfl{A}{} & &\left\|\vbox to
6truemm{}\right. &&\cr
0\rightarrow & P_0^*(t) &\hfl{M'}{} & P_0 &
\hfl{q}{} & {\cal F}
&\rightarrow 0&,\cr }$$where the vertical arrows are isomorphisms.
\ind We have $AM=M'B$, hence, since $M$ and $M'$ are $\varepsilon$\tx
symmetric, $M\,^t\!A={}^tBM'$, and therefore ${}^tBAM=M\,^t\!AB$. By
\ref{aut} this implies  ${}^t\!AB=\lambda I$ for some  $\lambda\in
k^*$. Multiplying
$A$ by a scalar we get $M'=AM\,^t\!A$. Thus all $\varepsilon$\tx
symmetric matrices providing a minimal resolution of ${\cal F}$ are
conjugate under the action of $\Aut(L_0)$. In the same way we see
that every automorphism of ${\cal F}$ is induced by a matrix
$A\in\Aut(L_0)$ such that $AM\,^t\!A=\lambda M$ for some $\lambda\in
k^*$. \medskip
\ind As above let us rephrase Theorem B in the way we will mostly use
it:
\th{Corollary} 
\enonce Assume $\ch(k)\not= 2$. Let $X$ be an integral
hypersurface of degree
$d$ in $\pn$,  and
$E$  an {\rm ACM} line bundle on $X$ with $E^2\cong {\cal
O}_X(d+t)$ {\rm (resp. }an {\rm ACM} rank $2$ vector bundle  on $X$
with determinant ${\cal O}_X(d+t)\,)$. There exists a symmetric {\rm
(resp. }skew-symmetric{\rm )} matrix
$M=(m_{ij})\in {\bf M}_\ell (\s)$, with $m_{ij}$ homogeneous of degree
$d_i+d_j-t$,  and an exact sequence
$$0\rightarrow  \sdir_{i=1}^{\ell }{\cal O}_{\pn}(t-d_i)\qfl{M} 
\sdir_{i=1}^{\ell }{\cal O}_{\pn}(d_i)\longrightarrow E\rightarrow
0\ ;$$ $X$ is defined by the equation $\det M=0$ {\rm (resp.} $\pf
M=0\,)$.
 If $H^0(X,E(-1))=0$ and $t=-1$, the matrix
$M$ is  linear, and the exact sequence takes the form
$$0\rightarrow {\cal O}_{\pn}(-1)^{rd}\qfl{M} {\cal
O}_{\pn}^{rd}\longrightarrow E\rightarrow 0$$with $r=\rk E$.
\endth\label{symdet}
\pr By assumption $E$ carries an $\varepsilon$\tx
symmetric form $E\otimes E\rightarrow {\cal O}_X(d+t)$, with
$\varepsilon=(-1)^{r-1}$. Then Theorem B provides the above minimal
resolution; by (\ref{hyper}) we have  $F=\det M$ if
$r=1$ and $F^2=\det M=(\pf M)^2$ if $r=2$. Moreover if $t=-1$ we have
$H^{n-1}(X,E(1-n))\cong H^0(X,E(-1))^*$ by Serre duality, so the last
assertion follows from Prop. \ref{linear}.\cqfd 

\section{Plane curves as determinants}
\ind  Let $C$ be a smooth plane curve of degree $d$, defined by an
equation $F=0$. We denote by $g={1\over 2}(d-1)(d-2)$ the genus of $C$.
Any line bundle $L$ on $C$ is ACM, hence admits a minimal
resolution (A1), with $\det M=F$.
\ind The case of  line bundles of degree $g-1$ follows directly from
Cor. \ref{lindet} (applied to $L(1)$):
\th Proposition
\enonce  {\rm a)} Let $L$ be a line bundle of degree $g-1$ on $C$ with 
$H^{0}(X,L)=0$. There exists a $d\times d$ linear
matrix $M$ such that $F=\det M$, and an exact
sequence
$$0\rightarrow {\cal O}_{\pl}(-2)^d\qfl{M} {\cal
O}_{\pl}(-1)^d\longrightarrow L\rightarrow 0\ .$$
\ind {\rm b)} Conversely, let $M$ be a $d\times d$ linear matrix such 
that $F=\det M$. Then the cokernel of $M: {\cal
O}_{\pl}(-2)^d\rightarrow {\cal O}_{\pl}(-1)^d$ is a line bundle $L$ on
$C$ of degree $g-1$ with $H^{0}(X,L)=0$.\cqfd
\endth\label{g-1}\smallskip 
\subsection Let  $|{\cal O}_{\pl}(d)|_{sm}$ be the open subset of the
projective space
$|{\cal O}_{\pl}(d)|$ parametrizing smooth plane curves of degree $d$.
For $\delta\in {\bf Z}$, let  ${\cal J}_d^{\delta}\rightarrow  |{\cal
O}_{\pl}(d)|_{sm}$ be the family of degree $\delta$ Jacobians: ${\cal
J}_d^\delta$ parametrizes pairs $(C,L)$ of a smooth plane curve of
degree $d$ and a line bundle of degree $\delta$ on $C$. Finally we
denote by
$\Theta_d$ the divisor in ${\cal J}_d^{g-1}$ consisting of pairs
$(C,L)$ with
$H^0(C,L)\not= 0$. It is an ample divisor, so its complement in 
${\cal J}_d^{g-1}$ is affine. 
\ind Let  ${\cal L}_d$ the open subset  of the vector
space of linear matrices  $M\in{\bf M}_d(\s^2)$ such that the equation 
$\det M=0$ defines a smooth plane curve  $C_M$ in $\pl$. By
associating to $M$  the curve $C_M$ and the line bundle
$L_M:=\Coker \bigl[{\cal O}_{\pn}(-2)^d\qfl{M} {\cal
O}_{\pn}(-1)^d\bigr]$ on
$C_M$ we define a morphism $\pi: {\cal L}_d\rightarrow {\cal
J}_d^{g-1}\moins\Theta _d$.  The group $GL(d)\times GL(d)$ acts on
${\cal L}_d$ by $(P,Q)\cdot M=PMQ^{-1} $; this action factors through
the quotient $G_d$ of $GL(d)\times GL(d)$ by 
${\bf G}_m$ embedded diagonally.
\th Proposition
\enonce The group $G_d$ acts freely and properly on ${\cal L}_d$;
the morphism
$\pi $ induces an isomorphism
${\cal L}_d/G_d\rightarrow {\cal
J}_d^{g-1}\moins\Theta _d$.
\endth\label{qtdet}
\pr  This is proved for instance in [B3], \S 3; let us give a proof
based on our present methods. Let
$M\in {\cal L}_d$,
$(P,Q)\in GL(d)\times GL(d)$, and 
$M'=PMQ^{-1}$. Then $\det M'=\det M$ up to a scalar, and we have a
commutative diagram
$$\diagram{0 &\longrightarrow &{\cal O}_{\pl}(-1)^d &\hfl{M}{} &{\cal
O}_{\pl}^d&\longrightarrow  &L_M&\rightarrow &0&\cr
 & &\vfl{Q}{} & & \vfl{P}{}&  &\vfl{}{\kern-4pt\wr}& &&\cr
0 &\longrightarrow &{\cal O}_{\pl}^d &\hfl{M'}{} &{\cal
O}_{\pl}(-1)^d&\longrightarrow  &L_{M'}&\rightarrow &0&;\cr
}\eqno(\ref{qtdet}.a)$$ 
thus $\pi$ factors through a morphism ${\cal L}_d/G_d\rightarrow 
{\cal J}_d^{g-1}\moins\Theta _d$. Conversely, if two matrices $M$ and
$M'$ give rise to isomorphic pairs, the minimal resolution of $L_M$ and
$L_M'$ are isomorphic, so we have a diagram $(\ref{qtdet}.a)$,
which shows that
$M$ and $M'$ are conjugate under $G_d$. Thus the orbits of $G_d$ in 
${\cal M}_d$ are isomorphic to the fibres of $\pi $, hence are closed.
Moreover by (\ref{aut}) the stabilizer of $M$ in $GL(d)\times GL(d)$
reduces to
${\bf G}_m$ embedded diagonally, hence $G_d$ acts freely on ${\cal
L}_d$.  This proves our assertions.\cqfd\smallskip

\rem{Remark}\label{birat} A simpler birational presentation of the
quotient
$GL(d) \backslash{\cal L}_d/GL(d)$ (and therefore of ${\cal
J}_d^{g-1}$) is obtained as follows. Let
${\cal D}_d$ be the closed subset of ${\cal L}_d$ consisting of
matrices of the form $X_0I_d+X_1M_1+X_2M_2$; it is isomorphic to an
affine open subset of
${\bf M}_d\times {\bf M}_d$, where ${\bf M}_d$ denotes the $k$\tx
variety of
$(d\times d)$\tx matrices. Then $G_d{\cal D}_d$ is an open
affine subset of  ${\cal L}_d$, and the stabilizer of
${\cal D}_d$ in $G_d$ is $PGL(d)$ acting on ${\bf M}_d\times {\bf
M}_d$ by conjugation. We thus have an open embedding ${\cal
D}_d/PGL(d)\mono GL(d) \backslash{\cal L}_d/GL(d)$.
\ind These quotients are of course unirational. It is a
classical question to decide whether they are rational: this would
have interesting applications in algebra (where the function field of
${\cal D}_d/PGL(d)$ is known as the ``center of the generic division
algebra") and in geometry (${\cal
D}_d/PGL(d)$ is birationally equivalent to the moduli space of stable
rank $d$ vector bundles on ${\bf P}^2$ with $c_1=0$, $c_2=d$). The
rationality is known only for
$d\le 4$. We refer to [L] for an excellent survey of these questions. 
\ind It is amusing to observe that the universal Jacobian ${\cal
J}_d^g$ {\it is rational} ([B3], 3.4): using the rational map
${\cal J}_d^g\dasharrow {\rm Sym}^g(\pl)$ which maps a general
pair $(C,L)$ to the unique element of $|L|$, we see that ${\cal
J}_d^g$ is birational to a projective fibre bundle over the
rational variety ${\rm Sym}^g(\pl)$. Unfortunately this does not seem
to have any implication on the more interesting question of the
rationality of
${\cal J}_d^{g-1}$.
\medskip
\ind We will now determine the minimal resolution  of a 
  generic line bundle $L$ of arbitrary degree on a generic
plane curve. Replacing $L$ by $L(t)$ for some $t\in{\bf Z}$ we
can assume $g-1\le
\deg L\le g-1+d$. 
\th Proposition
\enonce Let $L$ be a line bundle of degree $g-1+p$ on $C$, with $0\le
p\le d$. The following conditions are equivalent:
\indp{\rm (i)} $H^0(C,L(-1))=H^1(C,L)=0$, and the natural map
$$\mu_0:H^0(C,L)\otimes H^0(C,{\cal
O}_C(1))\rightarrow H^0(C,L(1))$$ is of maximal rank
{\rm (}that is, injective for
$p\le {d\over 2}$ and surjective for $p\ge {d\over 2}\,);$
\indp{\rm (ii)} There is an exact sequence
$$0\rightarrow {\cal O}_{\pl}(-2)^{d-p} \
\hfl{M}{}\ {\cal O}_{\pl}(-1)^{d-2p}\oplus {\cal
O}_{\pl}^{\,p}\longrightarrow 
 L\rightarrow 0\eqno{\hbox{if }\ p\le {d\over
2}\,,}$$ 
$$0\rightarrow {\cal O}_{\pl}(-2)^{d-p} \oplus {\cal
O}_{\pl}(-1)^{2p-d}\ \hfl{M}{}\ {\cal O}_{\pl}^{\,p}\longrightarrow 
 L\rightarrow 0\eqno{\hbox{if }\ p\ge {d\over 2}\,,}$$with
$\det M=F$.
\ind The set of pairs $(C,L)$ satisfying these conditions is Zariski
dense in 
${\cal J}_d^{g-1+p}$ {\rm (}and open if $k=\bar k\,)$.\endth
\label{gen} 
\pr  Assume that (i) holds. The natural maps
$$\mu_j:H^0(C,L(j))\otimes H^0(C,{\cal O}_C(1))\rightarrow
H^0(C,L(j+1))$$ are surjective for
$j\ge 1$ because
$H^1(C,L)=0$ [Md]; since $H^0(C,L(-1))=0$,
this means that the $\s^2$\tx module $\hbox{\ru G}_*(L)$ is generated
by homogeneous elements of degree 0 and 1. In other words,
 the minimal resolution of $L$ takes the
form $$0\rightarrow \sdir_{i=1}^{p+q}{\cal O}_{\pl}(e_i) \qfl{M} {\cal
O}_{\pl}(-1)^q\oplus{\cal O}_{\pl}^{\,p}\longrightarrow 
 L\rightarrow 0$$for some integer $q\ge 0$(observe that $\dim
H^0(C,L)=p$ by Riemann-Roch). The vanishing of $H^1(C,L)$ and the
minimality of the resolution imply 
$e_i\in\{-2,-1\}$, so we have
$$0\rightarrow  {\cal O}_{\pl}(-2)^{d-p}\oplus{\cal
O}_{\pl}(-1)^{\,r}\qfl{M} {\cal O}_{\pl}(-1)^q\oplus{\cal
O}_{\pl}^{\,p}
\longrightarrow L\rightarrow 0\eqno(\ref{gen}.a)$$
with $r=2p-d+q$. After tensor product with ${\cal O}_{\pl}(1)$
 the cohomology exact sequence gives $$  q=\dim \Coker \mu_0\ ,\ \
r=\dim \Ker
\mu_0\ ,\eqno(\ref{gen}.b)$$from which
(ii) follows.\smallskip 
\ind If (ii) holds,  we have the exact sequence
$(\ref{gen}.a)$ with $r=0$ (if $p\le {d\over 2}$) or $q=0$ (if $p\ge
{d\over 2}$). By $(\ref{gen}.b)$ $\mu_0$ is of maximal rank; the
vanishing of $H^0(C,L(-1))$ and $H^1(C,L)$  is clear.\smallskip 
\ind  Let $V$ be the  vector space of matrices $M$ appearing in (ii),
and $V_0$ the open subset of matrices whose determinant defines a smooth
curve; observe that $V_0$ is non-empty by (\ref{cursur}). As in
(\ref{qtdet}) we get a  morphism
$\pi :V_0\rightarrow {\cal J}_d^{g-1+p}$; since property (i) is open in 
${\cal J}_d^{g-1+p}$, $\pi $ is dominant. The last assertion of the
Proposition follows.\cqfd\medskip
\ind We just proved:
\th Corollary 
\enonce The variety ${\cal J}_d^\delta$ is unirational for all
$\delta\in{\bf Z}$.\cqfd
\endth\label{junirat}\smallskip 
\rem{Example} \label{d/2} Let us consider the relative Jacobian ${\cal
J}_d^0$.  We have $g-1=$ ${1\over 2}d(d-3)$, so if
$d$ is odd, ${\cal J}_d^0$ is canonically isomorphic to  ${\cal
J}_d^{g-1}$. Assume  $d=2e$, so that ${\cal J}_d^0$ is canonically
isomorphic to ${\cal J}_d^{g-1+{e}}$. For  $(C,L)$ generic in
${\cal J}_d^{g-1+{e}}$ 
 the minimal resolution of
$L$ takes the form
$$0\rightarrow {\cal O}_{\pl}(-2)^{e} \qfl{M}  {\cal
O}_{\pl}^{\,e}\longrightarrow 
 L\rightarrow 0\ ,$$so  the equation of $C$ can be written as 
the determinant of a matrix  $M\in{\bf M}_{e}(\s^2)$ with quadratic
entries. Writing such a matrix as
$M=\sum X_iX_j M_{ij}$, we see as in (\ref{birat})
 that  ${\cal J}_d^0$ is birationally equivalent to the quotient of
${\bf M}_{e}^5$ by 
$GL({e})$ acting by conjugation. This
quotient is birationally equivalent to a vector bundle over 
${\bf M}_{e}^2 /GL({e})$ [L]; in particular, we see
that the variety ${\cal J}_d^0$ is rational for $d=4,6$ or $8$.

\section{Plane curves as symmetric determinants}
\ind By Corollary \ref{symdet}, any line bundle $L$ on $C$ with
$L\ten\cong {\cal O}_C(s)$ admits a symmetric minimal resolution.
There are (at least) two interesting applications.\medskip

\subsec {\it Theta-characteristics}
\ind Recall that a {\it theta-characteristic} on a smooth curve $C$ is
a line bundle $\kappa$ such that $\kappa\ten\cong K_C$. We write
$h^0(\kappa):=\dim H^0(C,\kappa)$.

\th Proposition
\enonce Let $C$ be a smooth plane curve, defined by an equation
$F=0$, and $\kappa$ a theta-characteristic on $C$. 
\indp{\rm a)} If $h^0(\kappa)=0$, $\kappa$ admits a minimal resolution
$$0\rightarrow {\cal O}_{\pl}(-2)^d\qfl{M} {\cal
O}_{\pl}(-1)^d\longrightarrow \kappa\rightarrow 0\ ,$$where the matrix
$M\in{\bf M}_d(\s^2)$ is symmetric {\rm (}linear{\rm )}
and
$\det M=F$.
\indp{\rm b)} If $ h^0(\kappa)=1$, $\kappa$ admits a minimal resolution
$$ 0\rightarrow  {\cal O}_{\pl}(-2)^{d-3}\oplus
{\cal O}_{\pl}(-3)\ \hfl{M}{}\  {\cal
O}_{\pl}(-1)^{d-3}\oplus {\cal O}_{\pl}\longrightarrow
\kappa\rightarrow 0$$with a symmetric matrix $M\in{\bf
M}_{d-2}(\s^2)$ satisfying $\det M=F$,
 and  of  the form
$$M=\pmatrix{L_{1,1} & \cdots & L_{1,d-3} & Q_1\cr
\vdots & & \vdots & \vdots\cr
 L_{1,d-3} & \cdots & L_{d-3,d-3} & Q_{d-3}\cr 
Q_1 & \cdots & Q_{d-3} &H
}$$ 
where the forms $L_{ij}, Q_i, H$ are linear, quadratic and cubic
respectively.
\ind Conversely, the cokernel of a symmetric matrix $M$ as in {\rm a)}
{\rm (resp. b))} is a theta-characteristic $\kappa$ on $C$ with
$h^0(\kappa)=0$ {\rm (resp.} $h^0(\kappa)=1${\rm )}.
\endth\label{theta}
\ind Part a) is well-known, and goes back essentially to Dixon
[Di].  Part b) is stated for instance (without proof) in [B1],
6.27. Geometrically, when $\ch(k)\not= 2$, a) means 
that $C$ is the discriminant curve of a net of quadrics in ${\bf
P}^{d-1}$; b) means that $C$ is the discriminant curve of the quadric
bundle obtained by projecting the cubic hypersurface $\sum
U_iU_jL_{ij}+\sum U_iQ_i+H=0$ in the projective space ${\bf P}^{d-1}$
with coordinates
$U_1,\ldots,U_{d-3},X_0,X_1,X_2$ from the subspace $X_0=X_1=X_2=0$.
\smallskip  
\pr  Part a) follows directly from Cor. \ref{symdet} (applied
to $E=\kappa(1)$).
\ind Let $\kappa$ be a theta-characteristic on $C$, with
$h^0(\kappa)=1$. Then  $H^1(C,\kappa(1))=$ $H^0(C,\kappa(-1))^*=0$, so
$\hbox{\ru G}_*(\kappa)$ is generated by its elements of degree $0,1$
and $2$. In view of
\ref{symdet}, the minimal resolution of $\kappa$ is of the form
$$ 0\rightarrow {\cal O}_{\pl}(-1)^q\oplus {\cal O}_{\pl}(-2)^p\oplus
{\cal O}_{\pl}(-3)\ \hfl{M}{}\ {\cal O}_{\pl}(-2)^q\oplus {\cal
O}_{\pl}(-1)^p\oplus {\cal O}_{\pl}\longrightarrow \kappa\rightarrow
0$$
for some non-negative integers $p,q$. 
Since the resolution is minimal the summand ${\cal O}_{\pl}(-1)^q$ in
the first term is mapped into  ${\cal O}_{\pl}$; this implies $q\le
1$, and in fact $q=0$ because otherwise the non-zero section of
$\kappa$ would be annihilated by some linear form. This gives the form
of the resolution (and of the matrix $M$) in part b).\cqfd\medskip 
\ind Assume now $\ch(k)= 0$\note{This works equally well in all
characteristics $\not= 2$, but references are lacking.}. The moduli
space of pairs
$(C,\kappa)$, where
$C$ is a smooth plane curve of degree $d$ and $\kappa$ a
theta-characteristic on $C$, has two components, corresponding to the
parity of $h^0(\kappa)$, plus one special component when $d$ is odd
consisting of the pairs
$(C,{\cal O}_C((d-3)/2))$ ([B2], Prop. 3); a general element
$(C,\kappa)$ in a non-special component satisfies
$h^0(\kappa)\le 1$. 
\th Corollary
\enonce Each component of the  moduli space of  smooth plane curves
with a theta-characteristic is unirational.\cqfd
\endth\smallskip 
\rem{Remark} If $k$ is algebraically closed,  any smooth plane curve
admits a theta-char\-acteristic
$L$  with $H^0(L)=0$: this follows (via the Riemann singularity
theo\-rem) from the classical fact that the theta divisor of a
principally polarized Abelian variety cannot contain all points of
order 2 (see for instance [I], Ch. IV, lemma 11). Thus {\it every}
smooth plane curve can be defined by a symmetric linear determinant.
Actually every plane curve $C$ admits such a representation: one
reduces readily to the case when $C$ is integral; then one applies
Theorem B to the sheaf
$\pi_*L$, where $\pi: C'\rightarrow C$ is the normalization of $C$
and $L$ is a theta-characteristic on $C'$ with
$H^0(C',L)=0$\note{This remark answers a question of F.
Catanese.}.\medskip

\subsec {\it Half-periods}
\ind We assume again $\ch(k)=0\ ^2$. Let us consider
now the moduli space ${\cal R}_d$ of pairs
$(C,\alpha)$, where $C$ is a smooth plane curve of degree $d$ and
$\alpha$ a ``half-period", that is, a non-trivial line bundle on 
$C$ with $\alpha\ten\cong {\cal O}_C$. If $d$ is odd the map
$(C,\alpha)\mapsto (C,\alpha((d-3)/2))$ gives an isomorphism
of  ${\cal R}_d$ onto the above moduli space; we thus restrict
to  case $d$ even, say $d=2e$. It follows then from [B2], Prop.~2 
that ${\cal R}_d$ is {\it irreducible}.
\th Proposition
\enonce For $(C,\alpha)$ general in ${\cal R}_d$, the line
bundle $\alpha$ admits a minimal resolution
$$0\rightarrow {\cal O}_{\pl}(-e-1)^e \qfl{M}  {\cal
O}_{\pl}(-e+1)^e\longrightarrow 
 \alpha\rightarrow 0\ ,$$
where the matrix
$M\in{\bf M}_e(\s^2)$ is symmetric {\rm (}with quadratic
entries{\rm )} and $\det M=F$.
\endth
\pr In view of Cor. \ref{symdet}, this amounts to say
that the line bundle $\alpha(e-1)$ satisfies the equivalent
conditions of Prop. \ref{gen}. As in \ref{gen}, it suffices to exhibit
a symmetric matrix $M\in{\bf M}_e(\s^2)$ with quadratic
entries such that the equation $\det M=0$ defines a smooth plane curve.
\ind Start with a symmetric linear matrix $(L_{ij})\in {\bf M}_e(S)$  
such that the curve $\Gamma $ defined by $\det (L_{ij})=0$ is smooth
(such a matrix exists by Prop.
\ref{theta}). Changing coordinates if necessary we can assume that
$\Gamma $ is transverse to the coordinate axes and does not pass 
through the  intersection point of any two axes. Consider the covering
$\pi :\pl\rightarrow \pl$ given by $\pi
(X_0,X_1,X_2)=(X_0^2,X_1^2,X_2^2)$. The pull-back of $\Gamma$ by $\pi$
is smooth by our hypotheses; it is defined by the determinant of the
symmetric matrix $M= (L_{ij}(X_0^2,X_1^2,X_2^2))$ with quadratic
entries.\cqfd
\th Corollary
\enonce The moduli space ${\cal R}_d$ is unirational.\cqfd
\endth
\section{Plane curves as pfaffians}
\ind Again any rank $2$ vector bundle $E$ on the plane curve $C$ with
determinant ${\cal O}_C(s)$ for some integer $s$ admits a
skew-symmetric resolution. Let us restrict our attention to the
linear case.  Cor. \ref{symdet} applied to $E(1)$ gives:
\th Proposition
\enonce Let $C$ be a smooth plane curve of degree $d$, $E$ a rank $2$
vector bundle on $C$ with $\det E\cong K_C$ and $H^0(C,E)=0$. Then
$E$ admits a minimal resolution
$$0\rightarrow {\cal O}_{\pl}(-2)^{2d} \qfl{M} {\cal
O}_{\pl}(-1)^{2d}\longrightarrow E\rightarrow 0$$where the matrix
$M\in{\bf M}_{2d}(\s^2)$ is linear skew-symmetric  and $\pf M=F$.\cqfd
\endth
\ind Note that the condition $H^0(C,E)=0$ implies that $E$ is
semi-stable.  
\th Corollary
\enonce The moduli space of pairs
$(C,E)$, where $C$ is a smooth plane curve of degree $d$ and $E$ a
semi-stable rank $2$ vector bundle on $C$ with determinant
$K_C$, is  unirational.\cqfd
\endth
\ind This is not  surprising in this case, since the fibres
of the projection to $|{\cal O}_{\pl}(d)|$  are already
unirational.
\smallskip  
\subsection  Another consequence of the Proposition is that if
$d\ge 4$ and $M$ is general enough, the corresponding vector bundle
$E_M=\Coker M$ is stable and therefore simple, that is, $\End(M)=k$.
In view of (\ref{autsym}) this means that given 3 generic
skew-symmetric matrices $M_0,M_1,M_2\in {\bf M}_{2d}( k)$, the
equations
${}^t\!AM_iA=M_i$ for $i=0,1,2$ imply $A=\pm I$.
\label{simple}

\section{Surfaces as determinants} 
\subsection \label{AN}Let $S$ be a smooth surface of degree $d$ in
$\pt$, defined by an equation
$F=0$.  Let $C$ be a curve in $S$, and $L={\cal O}_S(C)$. Using the
exact sequence $0\rightarrow L^{-1} \rightarrow {\cal O}_S\rightarrow
{\cal O}_C\rightarrow 0$ and Serre duality, we see that $L$ {\it
is} ACM {\it if and only if $C$ is
  projectively normal in} $\pt$, that is, the restriction map
$H^0(\pt,{\cal O}_{\pt}(j))\rightarrow H^0(C,{\cal O}_{C}(j))$ is
surjective for every $j\in{\bf Z}$. Since any line bundle is
of the form  ${\cal O}_S(C)$ after some twist, this characterizes the
ACM line bundles on $S$. Thus any projectively normal curve
contained in $S$ gives rise to an expression of $F$ as the
determinant of a matrix $M\in{\bf M}_k(\s^3)$.  Recall however that
a  hypersurface section of $S$ gives the trivial case $M=(F)$;  a
curve $C$ defined in $\pt$ by two equations  $A=B=0$  produces 
a $2\times 2$\tx matrix $M=\pmatrix{A & B\cr C & D}$.\smallskip 
\ind We will now restrict our study to  {\it linear}
determinants.
\th Proposition
\enonce Let $C$ be a projectively normal curve on $S$, of degree
${1\over 2}d(d-1)$ and genus ${1\over 6}(d-2)(d-3)(2d+1)$. The line
bundle ${\cal O}_S(C)$ admits a minimal resolution
$$0\rightarrow {\cal O}_{\pt}(-1)^d\qfl{M} {\cal
O}_{\pt}^d\longrightarrow {\cal O}_S(C)\rightarrow 0$$with $\det M=F$.
\ind Conversely, let $M\in{\bf M}_d(\s^3)$ be a linear matrix such that
$\det M=F\,;$ the cokernel of $M:{\cal O}_{\pt}(-1)^d\longrightarrow 
{\cal O}_{\pt}^d$ is isomorphic to
${\cal O}_S(C)$, where $C$ is a smooth projectively normal curve on
$S$ with the above degree and genus.
\endth\label{surdet}
\pr  Let $C$ be a curve on $S$; put $L={\cal O}_S(C)$. A
straightforward Riemann-Roch computation   shows that the given
condition on the degree and genus of $C$ is equivalent to
$\chi(L(-1))=\chi(L(-2))=0$. If $C$ is projectively normal the
spaces $H^1(S,L(j))$ vanish (\ref{AN}), therefore the above condition
is also equivalent to $H^0(S,L(-1))=H^2(S,L(-2))=0$; this is exactly
what we need to apply Cor.~\ref{lindet}.
\ind Conversely, given a matrix $M$, let $L=\Coker M$; in
view of the above all we have to prove is that the linear
system $|L|$ contains a smooth curve. This is obvious in
characteristic $0$ since $L$ is spanned by its global
sections. In the general case, we first observe that the
restriction of $L$ to any smooth hyperplane section $H$ of $S$ is very
ample: indeed from the resolution $0\rightarrow {\cal
O}_{\pl}(-1)^d\rightarrow {\cal O}_{\pl}^d\rightarrow
L_{\,|H}\rightarrow 0$ we get $H^1(H,L_{\,|H}(-1))=0$, hence
$H^1(H,L_{\,|H}(-x-y))=0$ for all
$x,y\in H$. It follows that the linear system $|L|$ on $S$ separates
two points $x,y\in S$ (possibly infinitely close) unless the line 
$\langle x,y\rangle$ lies in $S$; in other words, the
morphism $\varphi_L:S\rightarrow {\bf P}^{d-1}$ defined by $|L|$
contracts finitely many lines, and  embeds the complement of these
lines. Then a general hyperplane  in ${\bf P}^{d-1}$ cuts down a
smooth curve $C\in|L|$.\cqfd\medskip

\subsection Under the hypotheses of the Proposition, Grothendieck
duality provides  a dual exact sequence (see the proof of Theorem B):  
$$0\rightarrow {\cal O}_{\pt}(-1)^d\qfl{{}^tM} {\cal
O}_{\pt}^d\longrightarrow L^{-1} (d-1)\rightarrow 0\ ;$$in other words,
{\it the involution $M\mapsto {}^tM$ on the space of linear matrices
corresponds to the involution $L\mapsto L^{-1} (d-1)$ on}
$\Pic(S)$.\label{dual}
\medskip
\ind As we already pointed out, a general form of degree $d$ on ${\bf
P}^3$ can be represented as a linear determinant only for $d\le 3$, the
only non-trivial case being $d=3$. There we find the following classical
result [G]:
\th Corollary
\enonce Assume $k$ is algebraically closed. A smooth cubic surface can
be defined by an equation
$\det M=0$, where $M$ is a $3\times 3$\tx linear matrix. There
are $72$ such representations {\rm (}up
to the action of $GL(3)\times GL(3)$ by left and right
multiplication{\rm ),} corresponding  in a one-to-one way to the
linear systems of twisted cubics on $S$.\cqfd   
\endth
\ind  There are various ways of describing the set of linear systems of
twisted cubics on $S$: they also correspond to the birational morphisms
$S\rightarrow \pl$, or to the sets of 6 lines on $S$ which do not
intersect each other. In terms of these, the involution $M\mapsto {}^tM$
corresponds to the Sch\"afli involution which associates to such a set
$\{\ell _1,\ldots,\ell _6\}$ the unique set $\{\ell '_1,\ldots,\ell
'_6\}$ such that the 12 lines $\ell _i,\ell '_j$ form a {\it double-six},
that is satisfy $\ell _i\cap \ell '_i=\varnothing$ and $\ell _i\cdot
\ell '_j=1$ for $i\not= j$. 

 \ind As a consequence, the space of pairs $(S,\lambda)$,
where $S$ is a smooth cubic surface and $\lambda$ a set of 6
  non-intersecting lines,  is {\it rational}: as in (\ref{birat}) it is
birational to the quotient of $({\bf M}_3)^3$ by the group $GL(3)$
acting by conjugation, and we know that this quotient is rational.
\smallskip 
\ind In the case of a non-necessarily algebraically closed field, we
find the following result of B. Segre [Se]:
\th Corollary
\enonce Let $S$ be a smooth cubic surface. The following conditions are
equivalent:
\indp{\rm (i)} $S$ can be defined by an equation $\det M=0$, where $M$
is a $3\times 3$\tx linear matrix;
\indp{\rm (ii)} $S$ contains a twisted cubic{\rm ;}
\indp{\rm (iii)} $S$ admits a birational morphism to $\pl\,;$
\indp{\rm (iv)} $S$ contains a rational point and a set {\rm (}defined
over
$k)$ of $6$  non-intersecting lines.
\endth
\pr The equivalence of (i), (ii) and (iii) follows from Prop.
\ref{surdet}. The implication (iii) $\Rightarrow $ (iv) is clear. If
(iv) holds, the surface obtained from $S$ by blowing down the set 
of $6$  non-intersecting lines is isomorphic to $\pl$ over $\bar k$ and
contains a rational point, hence is $k$\tx isomorphic to
$\pl$.\cqfd
\smallskip 
\th Corollary
\enonce A smooth quartic surface is determinantal if and only if it
contains a non-hyperelliptic curve of genus $3$, embedded in $\pt$ by
a linear system of degree $6$.
\endth
\pr  The only point to check is that a curve $C$ of genus $3$, 
embedded in $\pt$ by a linear system $|L|$ of degree $6$, is
projectively normal if and only if it is not hyperelliptic. Since
$H^1(C,L)=0$, the projective normality reduces using the base-point
free pencil trick to the surjectivity of the restriction map
$H^0(\pt,{\cal O}_{\pt}(2))\rightarrow H^0(C,L\ten)$; or equivalently,
since both spaces have the same dimension, to its injectivity. One
checks that $C$ is contained in a quadric if and only if it is
hyperelliptic.\cqfd\medskip
\subsection There is another approach to Prop. \ref{surdet}, which we
will now sketch. Given the linear matrix $M$, let $C$ be the
divisor of the section of $L=\Coker M$ corresponding to the first basis
vector of ${\cal O}_{\pt}^d$. Using (\ref{dual}) we see easily that
 the curve $C$ is defined in $\pt$ by the maximal minors of the
matrix $N$ obtained from $M$ by deleting the first row. Conversely,
since
 $C$ is projectively normal, it admits a resolution 
$$0\rightarrow \sdir_{j=1}^{\ell -1}{\cal
O}_{\pt}(e_j)\ \qfl{N}\ \sdir_{i=1}^{\ell} {\cal
O}_{\pt}(d_i)\,\qfl{\Delta}\, {\cal O}_{\pt}\rightarrow {\cal
O}_C\rightarrow 0\ ,$$
where $\Delta$ is given by the maximal minors of $N$; with some work
one finds $\ell =d$, $e_1=\ldots=e_{d -1}=-d$ and
$d_1=\ldots=d_{d}=-(d-1)$. It follows easily that any surface
of degree $d$  containing $C$ is defined by the determinant of a
linear matrix obtained by adding one  row to $N$.\label{minors}
\smallskip 
\subsection As indicated in the introduction, we will not
consider surfaces defined by symmetric determinants, though this is
again a classical and rich story; we refer to [C1] or [C2] for a
modern treatment.

\section{Surfaces as Pfaffians}
\ind\hskip1truecm {\it From now on we assume} $\ch(k)=0$ (see
\ref{charp}).
\subsection\label{spf} Again we will restrict ourselves to the linear
case, that is to surfaces $S\i\pt$ defined by an equation $\pf M=0$,
where $M$ is a $(2d)\times (2d)$ skew-symmetric linear matrix.
 
\ind Let  $Z$ be a finite reduced subscheme of $\pn$, of
degree\note{The degree of $Z$ is by definition $\dim_k H^0(Z,{\cal
O}_Z)$.}
$c$,
 and $\id_Z$ its homogeneous ideal in $\s^n$. 
 $Z$ is said to be {\it arithmetically Gorenstein} if the algebra
$\r:=\s/\id_Z$ is Gorenstein. This implies that there exists an integer 
$N$ such that:
\indp a) $\dim \r_p+\dim \r_{N-p}=c$ for all $p\in{\bf Z}$.
\ind The integer $N$ is uniquely determined: it is
the largest integer such that $\dim\r_N<c$.  By lack of a better name
we will call it the {\it index} of $Z$.  
\ind Assume $k=\bar k$. By [D-G-O], thm. 5, the subscheme $Z$ is
arithmetically Gorenstein if and only if it satisfies a) and:
\indp b) $ Z$ has the Cayley-Bacharach property w.r.t. the
linear system
$|{\cal O}_{\pn}(N)|$; that is, for each point $z\in  Z$,
every element of 
$|{\cal O}_{\pn}(N)|$ containing $ Z\moins z$ contains
$ Z$. 
\ind In general,  $Z$ is arithmetically Gorenstein if and
only if $Z\otimes_k\bar k$ has the same property.

\ind Let $Z\i\pt$ be a finite arithmetically Gorenstein subscheme, 
contained in a surface $S$ of degree $d$.  Let ${\cal I}_Z$ be the
sheaf of ideals of $Z$ in ${\cal O}_S$. Using the exact sequence
$0\rightarrow {\cal I}_Z\rightarrow {\cal O}_S\rightarrow {\cal
O}_Z\rightarrow 0$, property a) for $p=N$ gives
$\dim H^1(S,{\cal I}_Z(N))=1$. Thus there  exists a unique non-trivial 
extension (up to automorphism)
$$0\rightarrow {\cal O}_S\rightarrow   E\rightarrow {\cal
I}_Z(N-d+4)\rightarrow 0\ .$$
We claim that $E$ is locally free. To check this we can assume that $k$  
is algebraically closed;
then b) is equivalent to  $H^1(S,{\cal I}_{Z'}(N))=0$ for each proper 
subset $Z'\i  Z$, which implies our assertion by [G-H].
 We will say that $E$ is the vector bundle associated to $Z$.

\th Proposition
\enonce Let $S$ be a smooth surface of degree $d$ in $\pt$.
The following
conditions are equivalent:
\indp {\rm (i)} $S$ can be defined by an equation $\pf M=0$, where $M$   
is a skew-symmetric linear $(2d)\times (2d)$ matrix;
\indp {\rm (ii)} $S$ contains a finite arithmetically Gorenstein
reduced subscheme $Z$ of index $2d-5$, not contained in any surface of
degree $d-2$.
\ind More precisely, under hypothesis {\rm (ii),} the rank $2$ vector 
bundle $E$ associated to $Z$ admits a minimal resolution
$$0\rightarrow {\cal O}_{\pn}(-1)^{2d}\qfl{M} {\cal O}_{\pn}^{2d}
\longrightarrow E\rightarrow 0\ ;$$  
the degree of $Z$ is ${1\over 6}d(d-1)(2d-1)$.
\endth\label{surpf}
\pr If (i) holds the vector bundle $E:=\Coker M$ is spanned by its
global sections; let $Z$ be the zero locus of a general section of
$E$. Under (i) or (ii)  we have an exact sequence
$$0\rightarrow {\cal O}_S\rightarrow E\rightarrow {\cal
I}_Z(d-1)\rightarrow 0\ .\eqno(\ref{surpf}.a)$$ 
In view of Prop. \ref{symdet}, we have to prove the equivalence of:
\indp$\bullet$ $E$ is ACM and $H^0(S,E(-1))=0$;
\indp$\bullet$ $Z$ is arithmetically Gorenstein and 
$H^0(S,{\cal I}_Z(d-2))=0$.
\ind To do that we can assume $k=\bar k$. The fact that $E$ is  
locally free implies that
$Z$ has the  Cayley-Bacharach property w.r.t. 
$|{\cal O}_{\pt}(2d-5)|$ [G-H]. The sequence
$(\ref{surpf}.a)$ provides an isomorphism
$$H^0(S,E(-1))\iso H^0(S,{\cal I}_Z(d-2))\ ,$$  
and gives rise for each
$j\in{\bf Z}$ to an exact sequence
$$0\rightarrow H^1(S,E(j))\rightarrow  H^1(S,{\cal
I}_Z(d-1+j))\qfl{\partial} H^2(S,{\cal O}_S(j))\ .$$
Using the exact sequence $0\rightarrow {\cal I}_Z\rightarrow {\cal
O}_S\rightarrow {\cal O}_Z\rightarrow 0$, we can identify 
$H^1(S,{\cal I}_Z(k))$ with the cokernel of the restriction map
$r_k:H^0(S,{\cal O}_S(k))\rightarrow H^0(Z,{\cal O}_Z(k))$; the map 
$H^0(Z,{\cal O}_Z(d-1+j)\rightarrow H^2(S,{\cal O}_S(j))$ deduced
from $\partial$ is identified  by Serre duality to the transpose
of $r_{d-4-j}$. Therefore the vanishing of  $H^1(S,E(j))$ is
equivalent to $\im r_{d-1+j}=\Ker {}^tr_{d-4-j}= (\im
r_{d-4-j})^{\perp}$, that is to $\dim \r_{d-1+j}=c-\dim \r_{d-4-j}$.
This proves the equivalence of (i) and (ii).
\ind Under these equivalent conditions, we have $\Card Z=c_2(E)$; 
this number can be computed for instance using Riemann-Roch and
$\chi(E)=2d$.
\medskip
\rem{Remarks} a) We have to restrict to the  characteristic $0$ case 
because we do not know how to prove that the zero locus of a general
section of $E$ is smooth in characteristic $p$. The same problem 
occurs in higher dimension as well.\label{charp}
\ind b) As in (\ref{minors}) we could  use another approach: using the
 Buchsbaum-Eisenbud theorem [B-E] one shows that $\id_Z$ is
generated by
 the $(2d-2)\times (2d-2)$ pfaffians 
extracted from a skew-symmetric linear  $(2d-1)\times (2d-1)$\tx matrix
$N$; then $X$ is defined by the pfaffian of the matrix 
$\pmatrix{N & C \cr -^tC &0}$, where $C$ is a column of linear
forms.\smallskip 
\rem{Examples} For a cubic surface we have $\deg
Z=5$, and $N=1$. If $k=\bar k$ a subset $Z$ is arithmetically
Gorenstein if and only if any 4 points in $Z$ are linearly independent,
that is,
$Z$ is in general position.
\ind  For a quartic the subset $Z$ should have $14$ points,
not be contained in a quadric, and
satisfy the  Cayley-Bacharach property  w.r.t. $|{\cal
O}_S(3)|$.\label{expf}\smallskip 
\subsection Let us observe that for each  $d$ there exists smooth
surfaces defined by the pfaffian of a $(2d)\times (2d)$ skew-symmetric
linear matrix, and therefore containing a subset $Z$ with the
properties of the Proposition; we can take for instance 
$M=\pmatrix{0 & N\cr -{}^tN & 0}$, where $N$ is a linear $d\times d$
matrix: we have $\pf M=\det N$, and we can choose
$N$ so that the surface $\det N=0$ is smooth (\ref{cursur}). The
corresponding vector bundle $E$ is $L\oplus L^{-1} (d-1)$, where $L$ is
the line bundle $\Coker N$; the zero set $Z$  
 of a general section of $E$ is the intersection of two 
curves on $S$ of the type described in Prop. \ref{surdet} (see also
\ref{exist}).  
\smallskip 
\ind  We will now investigate when a generic surface of degree $d$
can be written as a linear pfaffian. 
\th Proposition
\enonce Assume that $k$ is algebraically closed. 
\ind {\rm a)} Every cubic surface can
be defined by a linear pfaffian.
\ind {\rm b)} A general surface of degree $d$ in $\pt$ can be
defined by a linear pfaffian if and only if $d\le 15$.
\endth\label{gensurpf}
\pr   a) follows from Proposition \ref{surpf} and Example
\ref{expf}.  Let ${\cal S}_d$ be the variety of linear 
skew-symmetric matrices  $M\in{\bf M}_{2d}(\s^3)$ such that  the
equation $\pf M=0$ defines a smooth surface  in $\pt$. Consider
the  map 
$\pf:{\cal S}_d\rightarrow   |{\cal O}_{\pt}(d)|$. 
We have $\dim{\cal S}_d/GL(2d) =
4d(2d-1)-4d^2= 4d(d-1)$; an easy computation gives $4d(d-1)<\dim 
|{\cal O}_{\pt(d)}|$ for $d\ge 16$, which gives the ``only if"
part of b). 
\ind To prove the remaining part we use Adler's method ([A-R], App.
V), namely we prove that the differential of $\pf$ is surjective at
a general matrix $M\in{\cal S}_d$. As in {\it loc. cit.}, a standard
computation shows that this is equivalent to the fact that the
vector space
$H^0(\pt,{\cal O}_{\pt}(d))$ is spanned by the forms
$X_k M_{ij}$, where $M_{ij}$ is
the pfaffian of the skew-symmetric matrix obtained from $M$ by
deleting the rows and columns of index $i$ and $j$. This has
been  checked by F.~Schreyer using the computer algebra
system Macaulay 2: a script is provided in the
Appendix.\cqfd\medskip
\ind We do not consider the proof of b) as completely satisfactory,
since it relies on a computer checking which does not  
provide any clue as why the result holds. The following lemma
explains better the meaning of the result. Recall that we associate
to a matrix $M\in{\cal S}_d$ the smooth surface
$S_M$ defined by $\pf M=0$ and the vector bundle $E_M:=\Coker \bigl[{\cal
O}_{\pn}(-1)^d\qfl{M} {\cal O}_{\pn}^d\bigr]$ on
$S_M$.
\th Lemma
\enonce The pfaffian map $\pf:{\cal S}_d\rightarrow   |{\cal O}_{\pt}(d)|$ 
is dominant if and only if $H^2(S_M,\en_0(E_M))$ vanishes for a general $M$ 
in ${\cal S}_d$.
\endth
\ind (As usual  $\en_0(E)$ denotes
the bundle of traceless endomorphisms of $E$.)
\pr We will restrict our attention to
matrices $M$ such that $E_M$ is  {\it simple}, that is, 
has only scalar endomorphisms. According to \ref{autsym}, this means
that  the only matrices  $A\in{\bf M}_d(k)$  such that $AM\,^t\!A=M$
are $\pm I$. The matrices $M$ with this property form an open
subset ${\cal S}_d^s$ of ${\cal S}_d$, which is non-empty by
\ref{simple}. 
\ind  We consider the map $\pf:{\cal S}_d^s\rightarrow  |{\cal
O}_{\pt(d)}|$; its fibre at a point $S\in |{\cal O}_{\pt(d)}|$ is
the moduli space of simple ACM rank 2 vector bundles on $S$ with
$\det E={\cal O}_S(d-1)$ and 
$H^{0}(S,E(-1))=0$.
A straightforward computation gives
$$\eqalign{\dim {\cal S}_d/GL(2d) &=\dim |{\cal O}_{\pt(d)}| 
-\chi(\en_0(E_M))\cr  &=\dim |{\cal O}_{\pt(d)}| +\dim
H^1(S_M,\en_0(E_M)) -\dim H^2(S_M,\en_0(E_M))}$$
\line{\hfill (7.7.{\it a})} for any
matrix
$M\in{\cal S}^s_d$.
\ind  If $H^2(S_M,\en_0(E_M))=0$, the moduli space of simple vector
bundles on $S_M$ is  smooth of dimension $\dim H^1(S_M,\en_0(E_M))$ 
at $[E_M]$. It then follows from (7.7.{\it a}) that
$\pf$ is dominant.
\ind Conversely, assume that $\pf$ is dominant. Let $S$ be a
generic surface of degree $d$; the
 fibre  $\pf^{-1} (S)$  can be identified with an open
subset of the moduli space of simple rank 2 bundles $E$ on $S$ with
$\det E={\cal O}_S(d-1)$ and $c_2(E)={1\over 6}d(d-1)(2d-1)$.
Being smooth,  this open subset is of dimension
$\dim H^1(S,\en_0(E))$. Comparing with (7.7.{\it a})
gives $H^2(S,\en_0(E))=0$.\cqfd\medskip

\subsection Thus assertion b) of Prop. \ref{gensurpf} is
equivalent to the fact that on a general surface $S$ of degree
$d$, the moduli space of  ACM rank 2 vector bundles  with $\det
E={\cal O}_S(d-1)$ and 
$H^{0}(S,E(-1))=0$ is smooth of the {\it expected dimension}
$-\chi(\en_0(E))$ for $d\le 15$. We were not able to prove
this directly, except for the obvious case $d=4$ where the
vanishing of $H^2(S,\en_0(E))$ follows from Serre duality.

\section{Threefolds  as  linear pfaffians}
\subsection Let us first briefly recall Serre's
construction. Let
$X$ be a projective manifold of dimension $\ge 3$ and $E$ a rank 2
vector bundle on $X$, spanned by its global sections; put $L=\det
E$.  Then the zero locus of a general section $s$ of $E$ is a  
 submanifold $V$ of codimension 2 in $X$, and there is an exact
sequence
$$0\rightarrow {\cal O}_X\qfl{s} E\longrightarrow {\cal
I}_VL\rightarrow 0\ ;$$
it follows that $K_V$ is isomorphic to $(K_X\otimes L)_{\,|V}$.
Conversely, given a codimension 2  submanifold $V\i X$ and a line
bundle $L$ on $X$ such that $K_V\cong (K_X\otimes L)_{\,|V}$,
there exists a rank 2 vector bundle $E$ and a section $s\in
H^0(X,E)$ such that
$Z(s)=V$; if moreover $V$ is connected, the pair $(E,s)$ is uniquely
determined up to isomorphism. We will refer to $E$ as the vector
bundle associated to $V$.
\ind Recall that a submanifold $V$ of $\pn$ is said to be {\it
arithmetically Cohen-Macaulay} if the sheaf ${\cal O}_V$ is ACM {\it
and} $V$ is projectively normal. This implies in particular
$H^0(V,{\cal O}_V)=k$, so $V$ is connected.

\th Proposition
\enonce Let $X$ be a smooth hypersurface of degree $d$ in $\pn$
$(n=4$ or$\, 5)$.  The following
conditions are equivalent:
\indp {\rm (i)} $X$ can be defined by an equation $\pf M=0$, where $M$   
is a skew-symmetric linear $(2d)\times (2d)$ matrix;
\indp {\rm (ii)} $X$ contains a codimension $2$ submanifold $V$ which
is arithmetically Cohen-Macaulay, not contained in any hypersurface of
degree $d-2$, and such that $K_V\cong {\cal O}_V(2d-2-n)$.
\ind More precisely, under hypothesis {\rm (ii),} the rank $2$
vector bundle $E$ associated to $V$ admits a minimal resolution
$$0\rightarrow {\cal O}_{\pn}(-1)^{2d}\qfl{M} {\cal O}_{\pn}^{2d}
\longrightarrow E\rightarrow 0\ ;$$  
the variety $V$ has degree ${1\over 6}d(d-1)(2d-1)$.
\endth
\pr If (i) holds the vector bundle $E:=\Coker M$ is spanned by its
global sections; let $V$ be the zero locus  of a general section of
$E$. Under (i) or (ii)  we have an exact sequence
$$0\rightarrow {\cal O}_X\rightarrow E\rightarrow
{\cal I}_V(d-1)\rightarrow 0\ .$$By Serre duality, $E$ is ACM if and
only if $H^i(X,E(j))=0$ for $1\le i\le n-3$; in view of the above
exact sequence this is equivalent to $V$ being arithmetically
Cohen-Macaulay.  Similarly the condition $H^0(X,E(-1))=0$
translates as $H^0(X,{\cal I}(d-2))=0$; we conclude by cor.
\ref{symdet}.
\ind The degree of $V$ can be deduced for instance from
(\ref{surpf}) by  restriction to a general 3-dimensional linear
subspace.\cqfd\label{n-folds}\medskip

\subsection Note that there exist indeed  smooth threefolds and 
fourfolds satisfying the equivalent conditions of the Proposition.
One way to see this is to consider the  vector space ${\bf
M}^{ss}_{2d}$ of skew-symmetric $(2d)\times (2d)$ matrices, and
the  universal pfaffian hypersurface ${\cal X}_d\i {\bf P}({\bf
M}^{ss}_{2d})$ consisting of singular matrices. The singular locus of
${\cal X}_d$ consists of those matrices which are of rank $\le 2d-4$,
and has codimension 6. Therefore for $n\le 5$  a generic $\pn\i{\bf
P}({\bf M}^{ss}_{2d})$ intersects ${\cal X}_d$ along a smooth
hypersurface in $\pn$, defined by the vanishing of a linear
pfaffian.\label{exist}
\medskip
\subsec {\it The cubic threefold}
\th Proposition
\enonce If $k=\bar k$, every smooth cubic threefold can be defined
by an equation
$\pf M=0$, where $M$ is a skew-symmetric linear $6\times 6$ matrix.
\endth
\ind As mentioned in the introduction, this result is due to Adler
([A-R], App. V) in the case of a {\it generic} cubic threefold.
\par
\pr Let $X$ be a smooth cubic threefold. In view of Prop. \ref{n-folds},
we have to prove that $X$ contains a normal elliptic quintic curve. 
This is essentially in [M-T], Remark 4.9; we sketch the argument since
the result we need is not explicitely stated there. We first observe
that
$X$ contains a non-normal elliptic quintic curve $C$ (that is,
contained in a hyperplane): in fact any smooth
hyperplane section $S$ of $X$ contains finitely many 5-dimensional
families of such curves (represent $S$ as ${\bf P}^2$ blown up at 6
points and consider the linear system of plane cubics passing through 4
of the 6 points). Varying the hyperplane section gives a 8-dimensional
family of non-normal elliptic quintic curves in $S$.
\ind Let $C$ be one of these curves; the normal bundle $N_{C/V}$ fits
into an exact sequence
$$0\rightarrow {\cal O}_C(1)\longrightarrow N_{C/V} \longrightarrow
N_{C/S}\rightarrow 0\ ,$$from which one deduces $H^1(C,N_{C/V})=0$ 
and $\dim H^0(C,N_{C/V})=10$. Therefore the Hilbert scheme of 
curves  of degree 5 and arithmetic genus $0$ in $V$ is smooth of
dimension 10 at
$C$. The general member of the component containing $C$ is a
smooth elliptic quintic  curve not contained in any hyperplane,
and therefore projectively normal.\cqfd\medskip

\subsection By Prop. \ref{symdet}, a rank 2 vector bundle $E$ on
$X$ is associated to a normal elliptic quintic if and only if 
$F=E(-1)$ satisfies $\det F={\cal O}_X$ and $H^0(X,F)=0$; since
$\Pic(X)={\bf Z}$, this last condition means that $F$ is {\it stable}
(with respect to  ${\cal O}_X(1)$).  Let ${\cal M}_X$ be the moduli
space of stable ACM rank 2 vector bundles on $X$ with trivial
determinant; it is smooth of dimension $5$ [M-T]. By a theorem of
Druel [Dr], this is also the moduli space
of stable rank 2 vector bundles on $X$ with $c_1=0$ and $c_2=2\ell $,
where $\ell $ denotes the class of a line in $H^4(X,{\bf Z})$; we will
not need this result here.
\ind Let us now vary $X$ and consider the space
${\cal M}$ of pairs $(X,F)$, where $X$ is a smooth element of 
$|{\cal O}_{{\bf P}^4}(3)|$ and $F\in {\cal M}_X$. By the
Proposition we have a dominant rational map from the space of
linear skew-symmetric matrices
$M\in {\bf M}_6(\s^4)$ onto the space ${\cal M}$, which is therefore
 {\it  unirational}.
\subsection Thanks to [M-T], this has the following nice consequence. We
now assume $k={\bf C}$. Let $|{\cal O}_{{\bf P}^4}(3)|_{sm}$ be the 
open subset of 
$|{\cal O}_{{\bf P}^4}(3)|$ parametrizing smooth cubic threefolds.
The intermediate Jacobians of cubic threefolds fit into a 
universal family
${\cal J}\rightarrow |{\cal O}_{{\bf P}^4}(3)|_{sm}$. More
generally, for each integer $k$ we can define a twisted
intermediate Jacobian $J^k(X)$, which parametrizes one-dimensional
cycles on $X$ with cohomology class
$k\ell$; this is a principal homogeneous space under the usual
intermediate Jacobian $J^0(X)$. These spaces fit into a family
${\cal J}^k$ over
$|{\cal O}_{{\bf P}^4}(3)|_{sm}$; while each $J^k(X)$ is isomorphic
to $J^0(X)$, it is not clear  that ${\cal J}^k$ is isomorphic to
${\cal J}$. However the class of a plane section is a canonical
element in each $J^3(X)$, giving a section of the fibration ${\cal
J}^3\rightarrow |{\cal O}_{{\bf P}^4}(3)|_{sm}$; this provides 
canonical  isomorphisms ${\cal J}^k\iso {\cal J}^{k+3}$ above
$|{\cal O}_{{\bf P}^4}(3)|_{sm}$. Note  also that for $p\in{\bf Z}$
the multiplication map ${\cal J}^k\qfl{\times p} {\cal J}^{pk}$ is
a  finite \'etale covering, since it is so on each fibre.
\th Corollary
\enonce The  intermediate Jacobian ${\cal J}$ of the universal family
of cubic threefolds is unirational.
\endth\label{ijunirat}
\pr  Associating to a pair $(X,F)$ in ${\cal M}$ the class of $c_2(F)$
defines a morphism ${\cal M}\rightarrow {\cal J}^2$ above $|{\cal
O}_{{\bf P}^4}(3)|_{sm}$. By [M-T] this morphism is \'etale, hence
dominant; thus ${\cal J}^2$ is unirational. Using the  maps ${\cal
J}^2\qfl{\times 3} {\cal J}^{6}\iso {\cal J}$, we conclude that
${\cal J}$ is unirational.\cqfd\medskip
\ind Let us discuss the case of higher degree threefolds. 
\th Proposition
\enonce Assume that $k$ is algebraically closed. 
 A general threefold of degree $d$ in ${\bf P}^4$ can be defined
by a linear pfaffian if and only if $d\le 5$.
\endth
\pr Let us denote  again by ${\cal S}_d$ the space of linear
skew-symmetric matrices $M\in {\bf M}_{2d}(\s^4)$ such that the
equation $\pf M=0$ defines a smooth hypersurface $X_M\i{\bf P}^4$.
As before the group $GL(2d)$ acts freely and properly on ${\cal
S}_d$, and the map $\pf:{\cal S}_d\rightarrow |{\cal O}_{{\bf
P}^4}(d)|$ factors through ${\cal S}_d/GL(2d)$.
\ind An easy computation gives $\dim {\cal S}_d/GL(2d)<\dim |{\cal
O}_{{\bf P}^4}(d)|$ for $d\ge 6$, so a general threefold of degree
$\ge 6$ is not pfaffian. For $d=4$ and $5$ one checks as in
\ref{gensurpf}  that the
differential of
$\pf$ at a generic matrix is surjective (Appendix; for $d=4$ this was
also observed in [I-M]).\cqfd\medskip

\subsection Exactly as in lemma 7.7 we find that the  map
$\pf:{\cal S}_d\rightarrow |{\cal O}_{{\bf P}^4}(d)|$ is dominant if
and only  if $H^2(X_M,\en_0(E_M))=0$ for $M$ general in ${\cal
S}_d$ -- that is, if the moduli space of the vector
bundles we are considering on a general quartic or quintic
threefold has the expected dimension. We see in particular that there
is a finite number of ways of representing a general
quintic as a pfaffian; this number is an instance of the
{\it generalized Casson invariant} considered by Thomas [T]. It
would be of course quite interesting to determine it. 
\section{Fourfolds as linear pfaffians}
\subsection Let us keep  the  notation of (\ref{3f}) for {\it
fourfolds} in ${\bf P}^5$. We find in this case  $\dim{\cal
S}_d/GL(2d)< \dim |{\cal O}_{{\bf P}^5}(d)|$ for $d\ge 3$, so  a
general hypersurface of degree 
$\ge 3$ in ${\bf P}^5$ cannot be defined by the vanishing of a linear
pfaffian  (a smooth hyperquadric can of course, since it is
isomorphic to the Grassmannian of lines in $\pt$ in the Pl\"ucker
embedding). For $d=3$ one finds
$\dim {\cal S}_3/GL(6)=\dim |{\cal O}_{{\bf P}^5}(3)|-1$.
\th Proposition
\enonce {\rm a)} A (smooth) cubic fourfold $X\i{\bf P}^5$ is
pfaffian if and only if
it contains a Del Pezzo surface of degree $5$.
\ind {\rm b)} Assume $k={\bf C}$. The map $\pf:{\cal
S}_3/GL(6)\rightarrow |{\cal O}_{{\bf P}^5}(3)|$ is generically
injective. In particular, pfaffian cubic fourfolds form a
hypersurface in  the space of all smooth cubic fourfolds. 
\endth
\def\alt{\mathop{\rm Alt}\nolimits}
\ind The pfaffian cubics play a key role in the proof that the 
variety of lines  contained in a cubic fourfold is irreducible
symplectic [B-D]. Cubic fourfolds containing a Del Pezzo surface of
degree $5$ were already considered by Fano [F].\smallskip 
\pr Part a) follows at once from Prop. \ref{n-folds}; let us prove
part b). 
\ind Let us introduce a 6-dimensional vector space $V$ and the
space $\alt (V)$ of bilinear alternate forms on $V$; we will
 view ${\cal S}_3$ as an open subset of
$\alt(V)^6=$ $\alt(V)\otimes_kk^6$. The map $\pf:{\cal S}_3\rightarrow
|{\cal O}_{{\bf P}^5}(3)|$ associates to  a sextuple $(\varphi
_0,\ldots,\varphi _5)$ the hypersurface $\pf(\sum_iX_i\varphi _i)=0$. 
The group $GL(6)$ acts on ${\cal S}_3$  through its  action on $k^6$;
this action commutes with the action of $GL(V)$, and  the map
$\pf:{\cal S}_3/GL(V)\rightarrow |{\cal O}_{{\bf P}^5}(3)|$ is
$GL(6)$\tx equivariant. The orbits of
$GL(6)$ in ${\cal S}_3$ correspond to 6-dimensional subspaces
$L\i\alt(V)$; to such a subspace is associated the isomorphism class
of the cubic hypersurface $X_L$ of degenerate forms in ${\bf P}(L)$. 
Since the action of $GL(6)$  is generically free on $|{\cal O}_{{\bf
P}^5}(3)|$, it
 is  sufficient to prove that the isomorphism class of $X_L$ 
determines  $L$ (up to the action of $GL(V)$).
\ind The orthogonal
$L^{\perp}$ of $L$ in $\ext^2V$ is 9-dimensional; the
locus of rank 2 bivectors in ${\bf P}(L^{\perp})$ is a K3 surface $S$ 
of genus 8 [B-D]. By [M], a general K3 surface of genus 8 is
obtained in this way, and this representation is unique: the
surface $S$ determines the space $L^{\perp}\i \ext^2V$ (and
therefore also the space $L\i
\alt(V)$) up to the action of $GL(V)$. So what we need to prove
is that {\it the cubic $X_L$ determines the $K3$ surface} $S$  up
to projective isomorphism.

\ind We proved in [B-D] that the variety $F$ of lines contained
in $X_L$ is a (complex) symplectic manifold, isomorphic to the
Hilbert scheme $S^{[2]}$; in particular, the group $H^2(F,{\bf Z})$
carries a canonical quadratic form, and there is a Hodge isometry 
$$H^2(F,{\bf Z})\iso H^2(S,{\bf Z})\sdir_{}^{\perp} {\bf Z}\delta\ 
,$$ where $H^2(S,{\bf Z})$ is endowed with the
intersection form and $\delta$ is a class of type $(1,1)$ and square 
$-2$. The polarization of $F$ given by the embedding in the
Grassmannian ${\bf G}(2,6)$ corresponds to the class
$2l-5\delta$, where $l$ is the polarization on $S$ deduced from the
embedding $S\i{\bf P}(L^{\perp})$.
\ind Let $L$ and $L'$ be two subspaces  of $\alt(V)$ which produce
isomorphic cubics; let $(S,l)$ and $(S',l')$ be the corresponding
polarized K3 surfaces. We then have a Hodge isometry
$$\varphi : H^2(S,{\bf Z})\sdir_{}^{\perp} {\bf Z}\delta \iso
H^2(S',{\bf Z})\sdir_{}^{\perp} {\bf Z}\delta'$$
which maps the class $2l-5\delta$ to the corresponding class
$2l'-5\delta'$. Assume $\Pic(S)={\bf Z}l$. Then we have $\Pic(S')=
{\bf Z}l'$, and $\varphi$ induces an isometry ${\bf Z}l\oplus {\bf
Z}\delta\iso$ ${\bf Z}l'\oplus {\bf Z}\delta'$ which maps
$2l-5\delta$ onto
$2l'-5\delta'$; an easy computation shows that this implies $\varphi
(\delta)=\varphi (\delta')$. Thus $\varphi $ induces a Hodge
isometry of $H^2(S,{\bf Z})$ onto $H^2(S',{\bf Z})$ mapping $l$ to 
$l'$;  by the Torelli theorem for K3 surfaces this implies that $(S,l)$
and $(S',l')$ are isomorphic.\cqfd
\vfill\eject
\centerline{\bf Appendix}\smallskip \centerline{\it Hypersurfaces
are generically pfaffian in the expected range} 
\centerline{Frank-Olaf Schreyer\note{Fakult\"at f\"ur Mathematik 
und Physik, Universit\"at Bayreuth,
D-95440 Bayreuth, Germany;\hfill\break
\hbox{\it E-mail
address} : {schreyer@btm8x5.mat.uni-bayreuth.de}}}
\bigskip
\ind In this appendix we prove  by a {\it Macaulay
2} computation that a generic surface of degree $d\le 15$ in $\pt$,
and a general threefold of degree $d\le 5$ in ${\bf P}^4$, can be
defined by the pfaffian of a skew-symmetric $2d\times 2d$ matrix
with linear entries (Propositions \ref{gensurpf} and  \ref{3f} in
the text). As explained in the paper, it is sufficient to prove
that for some matrix $M$ of this type the space of homogeneous
forms of degree $d$ is equal to ${\goth m}\cdot  \hbox{\tt
pfaffians}(2d-2,M)$, where ${\goth m}$ is the ideal
spanned by the coordinates and $\hbox{\tt pfaffians}(2d-2, M)$
 the ideal of submaximal pfaffians of $M$. We compute
the dimension of the latter space at randomly choosen skew
symmetric matrices over a finite field using {\it Macaulay 2}
[G-S].  The computation is within the range of nowadays
computers. On the computer ``alice'' of the Mathematical
Science Research Institute at Berkeley the following code was
executed in about 2 hours of cpu. The output verifies the
result.\par

{\baselineskip=13pt\parindent=12pt
\let\\=\par
\bigskip
\noindent
{\tt
isPrime(31991) \\
kk=ZZ/31991 -- this is a field \\
\bigskip\nospacedmath
\noindent
randomSkewMatrix = (e,S) -> ( \\
\indent -- returns a random e x e  skew symmetric matrix \\
\indent -- with linear entries in the ring S \\
\indent N:=binomial(e,2); \\
\indent R:=kk[t\_0..t\_(N-1)];\\
\indent G:=genericSkewMatrix(R,t\_0,e); \\
\indent substitute(G,random(S$\kern1pt {\bf
\hat{}}\,\{$0$\}$,S$\kern1pt {\bf
\hat{}}\,\{$N:-1$\}$)) \\
\indent ) -- end randomSkewMatrix \\
\bigskip
\noindent
subPfaffiansViaSyzygies = (M) -> ( \\
\indent -- This is an alternative to the command pfaffians(2d-2,M). \\
\indent -- It returns the generators of the ideal of the 2d-2 pfaffians  \\
\indent -- of the linear 2d x 2d skew symmetric matrix M  computed \\
\indent -- using the structure theorem of [B-E]: \\
\indent -- Under a mild genericity condition on the submatrix M1 \\
\indent -- the syzygies of the 2d-1 x 2d-1 skew matrix M1 are its
2d-1 \\
\indent -- principal pfaffians. \\
\indent -- If the computation fails, then the standard way is used. \\
\indent d:=lift((rank source M)/2,ZZ); \\
\indent syzygiesGivePfaffians=true; i:=0; S:=ring M; \\
\indent J:=generators ideal{0\_S}; \\
\indent while syzygiesGivePfaffians==true and (i<(2*d)) do (\\
\indent\indent      -- take i-th 2d-1 x 2d-1 skew submatrix     \\
\indent\indent     
M1:=transpose((transpose(M\_$\{$0..(i-1),(i+1)..(2*d-1)$\}$))
\\
\indent\indent      \_$\{$0..(i-1),(i+1)..(2*d-1)$\}$); \\
\indent\indent      N1:=syz(M1,DegreeLimit=>d); \\
\indent\indent      syzygiesGivePfaffians=((degrees source N1) ==
$\{\{$d$\}\}$); \\
\indent\indent      if syzygiesGivePfaffians==true then \\
\indent\indent          J=(J|flatten(N1)); \\
\indent\indent      i=i+1; \\
\indent ); \\
\indent if syzygiesGivePfaffians then (mingens image J) \\
\indent else (mingens image pfaffians(2*d-1,M)) \\
\indent ) -- end subPfaffiansViaSyzygies \\
\bigskip
\noindent
isDominant=(r,d) -> ( \\
\indent       S:=kk[x\_0..x\_r];  M:=randomSkewMatrix(2*d,S); \\
\indent       J:=subPfaffiansViaSyzygies(M); \\
\indent       N=syz(J,DegreeLimit=>d); \\
\indent       -- DegreeLimit=> d is carefully choosen to compute only\\
\indent       -- linear sysygies. From this the number of kk-linear\\
\indent       ---independent elements of degree d in the ideal\\
\indent       -- with generated by J can be computed:\\
\indent       cd=binomial(d+r,r)-(r+1)*rank(target N)+(rank source
N); \\
\indent       cd==0) -- end isDominant \\
\bigskip
\noindent
lowerBoundForDominantDegree = (r) -> ( \\
\indent dominant:=true; d:=2; \\
\indent while dominant do  \\
\indent     (d=d+1;dominant=isDominant(r,d);); \\
\indent d-1) \\
\bigskip
\noindent
isDominant(5,3) \\
cd \\
time d4=lowerBoundForDominantDegree(4) \\
time d3=lowerBoundForDominantDegree(3) \\
}\par}
\bigskip

\ind Note that we used the method to compute pfaffians via
syzygies, since this is faster than the command {\tt
pfaffians(2*d-2,M)}. The reason is that syzygy computations are
fast, while the {\tt pfaffian} command does not utilize much
special structure. For comments on the commands and the
{\it Macaulay 2} language we refer to the on-line help.
\ind Notice that the computation also shows that the closure of the
scheme of pfaffian cubic 4-folds form a hypersurface in $|{\cal
O}_{{\bf P}^5}(3) |$.

\vfill\eject
\centerline{REFERENCES} \vglue15pt\baselineskip12.8pt
\def\num#1{\smallskip\item{\hbox to\parindent{\enskip [#1]\hfill}}}
\parindent=1.4cm 
\num{A-R} A. {\pc ADLER}, S. {\pc RAMANAN}: {\sl Moduli of Abelian
varieties}. Lecture Notes in Math. {\bf 1644}, Springer-Verlag (1996).
\num{B1} A. {\pc BEAUVILLE}: {\sl 	Vari\'et\'es de Prym et jacobiennes
interm\'ediaires.} Ann. Sci. \'Ec. Norm. Sup. {\bf 10}, 309--391
(1977). 
\num{B2} A. {\pc BEAUVILLE}: {\sl 	Le groupe de monodromie des 
familles universelles d'hyper\-surfaces et d'inter\-sec\-tions
compl\`etes.} Complex analysis and Algebraic Geometry, LN {\bf 1194},
195--207;  Springer-Verlag (1986).
\num{B3} A. {\pc BEAUVILLE}: {\sl Jacobiennes des courbes spectrales 
et syst\`emes hamiltoniens compl\`etement int\'e\-grables.}  Acta
math. {\bf 164}, 211--235 (1990).
\num{B-D} A. {\pc BEAUVILLE}, R. {\pc DONAGI}: {\sl La vari\'et\'e des
droites d'une hypersurface cubique de dimension $4$}. 
C.R. Acad. Sc. Paris {\bf 301}, 703--706 (1985). 
\num{B-E} D. {\pc BUCHSBAUM}, D. {\pc EISENBUD}:
{\sl Algebra structure on finite free resolutions and some
structure theorem for ideals of codimension $3$}. Am. J. Math. {\bf
99}, 447--485  (1977).

\num{Bo} N. {\pc BOURBAKI}: Alg\`ebre, Ch. {\bf 10} (alg\`ebre
homologique). Masson (1980).
\num{C1} F. {\pc CATANESE}: {\sl Babbage's conjecture, contact of 
surfaces, symmetric determinantal varieties and applications}. Invent.
Math. {\bf 63}, 433--465  (1981). 
\num{C2} F. {\pc CATANESE}: {\sl Homological algebra and algebraic
surfaces}. Algebraic geo\-metry (Santa Cruz 1995), 3--56, Proc.
Symp. Pure Math. {\bf 62} I,  AMS (1997). 
\num{Ca} A. {\pc CAYLEY}: {\sl A memoir on quartic surfaces}.
Proc. London Math. Soc. {\bf 3}, 19--69 (1869-71).
\num{C-T} R.\ {\pc COOK}, A.\ {\pc THOMAS}: {\sl Line bundles and
homogeneous matrices}. Quart. J. Math. Oxford  {\bf 30}, 423--429 
(1979). 

\num{D-G-O} E. {\pc DAVIS}, A. {\pc GERAMITA}, F. {\pc ORECCHIA}:
{\sl  Gorenstein algebras and the Cayley-Bacharach theorem}.
Proc. Amer. Math. Soc.  {\bf 93},  593--597  (1985). 
\num{D} L.E. {\pc DICKSON}: {\sl Determination of all general
homogeneous polynomials expressible as  determinants with linear 
elements}. Trans. Amer. Math. Soc. {\bf 22}, 167--179 (1921).
\num{Di} A.C. {\pc DIXON}: {\sl Note on the reduction of a ternary
quantic to a symmetric determinant}. Proc. Cambridge Phil. Soc. {\bf
11}, 350--351 (1902).
\num{Dr} S. {\pc DRUEL}: {\sl Espace des modules des faisceaux 
semi-stables de rang $2$ et de classes de Chern $c_{1}=0$, $c_{2}=2$
et $c_{3}=0$ sur une hypersurface cubique lisse de}
${\bf P}^{4}$. Preprint math.AG/0002058.
\num{F} G. {\pc FANO}: {\sl Sulle forme cubiche dello spazio a cinque
dimensioni contenenti rigate razionale del $4^o$ ordine}. Comment. 
Math. Helvetici {\bf 15}, 71--80 \hbox{(1942-43)}.
\num{G} H. {\pc GRASSMANN}: {\sl Die stereometrischen Gleichungen
dritten grades, und die dadurch erzeugten Oberfl\"achen}. J.
Reine Angew. Math.  {\bf 49},  47--65 (1855).
\num{G-H} P. {\pc GRIFFITHS}, J. {\pc HARRIS}: {\sl Residues and 
zero-cycles on algebraic varieties}. Ann. of Math. {\bf 108},
 461--505  (1978). 
\num{G-S} D. {\pc GRAYSON},  M. {\pc STILLMAN}: {\sl Macaulay 
$2$},  {\tt http://www.math.uiuc.edu/ Macaulay2/}.
\num{H} O. {\pc HESSE}: {\sl Ueber determinanten und ihre
Anwendung in der Geometrie, insbesondere auf Curven vierter
Ordnung}. J. Reine Angew. Math. {\bf 49},  243--264 (1855).
\num{I} J.-I. {\pc IGUSA}: {\sl Theta functions}. Grundlehren der
math. Wiss. {\bf  194}. Springer-Verlag (1972). 
\num{I-M} A. {\pc ILIEV}, D. {\pc MARKUSHEVICH}: {\sl Quartic 3-fold:
Pfaffians, instantons and half-canonical curves}. Preprint
math/9910133.
\num{L} L. {\pc LE} {\pc BRUYN}: {\sl Centers of generic division 
algebras, the rationality problem 1965--1990}. Israel J. Math. {\bf
76},  97--111  (1991). 
\num{M} S. {\pc MUKAI}: {\sl Curves, $K3$ surfaces and Fano
$3$-folds of genus} 
$\leq 10$. Algebraic geometry and commutative algebra (in honor of M. 
Nagata), Vol. I, 357--377, Kinokuniya, Tokyo (1988). 
\num{Md} D. {\pc MUMFORD}: {\sl Lectures on curves on an
algebraic  surface}.  Annals of Math.
Studies {\bf 59}. Princeton University Press, Princeton (1966).
\num{M-T} D. {\pc MARKUSHEVICH}, A. {\pc TIKHOMIROV}: 
{\sl The Abel-Jacobi map of a moduli component of vector bundles
on the cubic threefold}. Preprint 
math.AG/ 9809140, to appear in
J. of Alg. Geom.
\num{S} F. {\pc SCHUR}: {\sl Ueber die durch collineare
Grundgebilde erzeugten Curven und Fl\"achen}. Math. Annalen {\bf
18}, 1--32 (1881).
\num{Se} B. {\pc SEGRE}: {\sl On the rational solutions of
homogeneous cubic equations in four variables}. Math. Notae
{\bf 11}, 1--68 (1951). 
\num{T} R.P. {\pc THOMAS}: {\sl A holomorphic Casson invariant for 
Calabi-Yau 3-folds, and bundles on $K3$ fibrations}. Preprint
math/9806111.
\vskip1cm
\def\pc#1{\eightrm#1\sixrm}
\hfill\vtop{\eightrm\hbox to 5cm{\hfill Arnaud {\pc BEAUVILLE}\hfill}
 \hbox to 5cm{\hfill DMA -- \'Ecole Normale
Sup\'erieure\hfill} \hbox to 5cm{\hfill (UMR 8553 du CNRS)\hfill}
\hbox to 5cm{\hfill  45 rue d'Ulm\hfill}
\hbox to 5cm{\hfill F-75230 {\pc PARIS} Cedex 05\hfill}}
\end